\documentclass[10pt,twoside]{amsart}
\usepackage{amsmath}
\usepackage{amssymb}
\usepackage{mathrsfs}		
\usepackage[curve]{xypic}

\usepackage[colorlinks=true, urlcolor=blue,bookmarks=true,bookmarksopen=true,
citecolor=blue,hypertex]{hyperref}

\numberwithin{equation}{section}

\newtheorem{theorem}[subsection]{Theorem}
\newtheorem{defn}[subsection]{Definition}
\newtheorem{assumption}[subsection]{Assumption}

\newtheorem{corollary}[subsection]{Corollary}

\newtheorem{convention}[subsection]{Convention}

\newtheorem{lemma}[subsection]{Lemma}
\newtheorem{prop}[subsection]{Proposition}
\newtheorem{remark}[subsection]{Remark}

\def \begineq{\begin{equation}}
\def \endeq{\end{equation}}

\def \bb{\mathbb}

\def \mc{\mathcal}
\def \mf{\mathfrak}
\def \ms{\mathscr}

\newcommand\hmf[1]{{\hat{\mf #1}}}
\newcommand\tmf[1]{{\tilde{\mf #1}}}

\def \CC{{\bb{C}}}

\def \GG{{\bb G}}
\def \JJ{{\bb{J}}}

\def \RR{{\bb{R}}}
\def \TT{{\bb{T}}}
\def \UU{{\bb{U}}}
\def \VV{{\bb{V}}}
\def \WW{{\bb{W}}}
\def \ZZ{{\bb{Z}}}

\def \GGC{{\mc G}}
\def \JJC{{\mc J}}
\def \KKC{{\mc K}}
\def \LLC{{\mc L}}

\def \NNC{{\mc N}}

\def \TTC{{\mc T}}
\def \XXC{{\mc X}}

\def \GGS{{\ms G}}
\def \XXS{{\ms X}}

\def \({\left(}
\def \){\right)}
\def \<{\langle}
\def \>{\rangle}
\def \ldb{{[\![}}
\def \rdb{{]\!]}}

\def \bar{\overline}

\def \dsum{\oplus}
\def \dto{\dashrightarrow}

\def \inter{\cap}
\def \into{\hookrightarrow}

\def \roof{{\hat{}}}
\def \tar{\underline}
\def \tensor{\otimes}
\def \union{\cup}

\def \xto{\xrightarrow}

\def \ad{{\rm ad}}
\def \Ann{{\rm Ann}}

\def \Diff{{\rm Diff}}

\renewcommand{\1}{1\!\!1}

\def \qed{\hfill $\square$ \vspace{0.03in}}

\begin{document}

\title{Reduction and duality in generalized geometry}
\author{Shengda Hu}
\email{shengda@dms.umontreal.ca}
\address{D\'epartement de Math\'ematiques et de Statistique, 
Universit\'e de  Montr\'eal, CP 6128 succ Centre-Ville, Montr\'eal, QC H3C 3J7,
Canada}

\abstract
Extending our reduction construction in \cite{Hu} to the Hamiltonian action of a
Poisson Lie group, we show that generalized K\"ahler reduction exists even when
only one generalized complex structure in the pair is preserved by the group
action. We show that the constructions in string theory of the (geometrical)
$T$-duality with $H$-fluxes for principle bundles naturally arise as reductions
of factorizable Poisson Lie group actions. In particular, the group may be
non-abelian.
\endabstract
\maketitle


\section{Introduction} \label{intro}

In this article, we propose a candidate of geometric realization of part of the
ansatz of $T$-duality with $H$-flux in the physics literature, using reductions
in generalized K\"ahler geometry. 
$T$-duality has long been intensively studied in physics and has made its marks
in mathematics as well, e.g. via mirror symmetry \cite{Strominger}. 
The context of our reduction construction is the Hamiltonian Poisson action of
Poisson Lie group. Classically, such reduction in symplectic category was first
discussed in \cite{Lu} and our construction here should be viewed as the
generalization of it to generalized geometry. 

Generalized geometry is introduced by Hitchin \cite{Hitchin} in the context of
generalized Calabi-Yau manifolds. The general theory of generalized complex and
K\"ahler geometries is first developed by Gualtieri in his thesis
\cite{Gualtieri}. Various reduction constructions in the context of generalized
geometry are developed by \cite{Bursztyn, Hu, Lin, Stienon, Vaisman}. The
approach taken here follows the point of view of Hamiltonian symmetries \cite{Hu}.

It is by now well-known that a generalized complex structure induces a canonical
Poisson structure, e.g.  \cite{Abouzaid, Crainic, Gualtieri, Hu}, also
\S\ref{poisson:induce}. Let $G$ be a Poisson Lie group with dual group $\hat
G$, then the Hamiltonian Poisson action with moment map as defined in \cite{Lu}
(also see definition \ref{app:poissonmom}) can be adapted to generalized complex
geometry (definition \ref{poisson:ham}), as well as generalized K\"ahler
geometry (definition \ref{kahler:ham}). We then have the first results on
reduction:

\vspace{0.1in}
\noindent
{\bf Theorem \ref{poisson:thm}, \ref{kahler:reduct}.} {\it Suppose $(M, \JJC)$
is an extended complex manifold with Hamiltonian $G$-action, whose moment map is
$\mu : M \to \hat G$. Let $M_0 = \mu^{-1}(\hat e)$, where $\hat e \in \hat G$ is
the identity element. Suppose that $\hat e$ is a regular value and the
geometrical action of $G$ is proper and free on $M_0$. Then there is a natural
extended complex structure on the reduced space $Q = M_0/G$.

If furthermore, $(M, \JJC_1, \JJC_2)$ is an extended K\"ahler manifold and the
$G$-action is $\JJC_1$-Hamiltonian. Then there is a natural extended K\"ahler
structure on the reduced space $Q$.
}

\vspace{0.1in}
\noindent
The notion \emph{extended} ($+$ \emph{structures}) is adopted to emphasize that
we consider $\TTC M$ as an extension of $TM$ by $T^*M$, instead of as a direct
sum, with an exact Courant algebroid structure (cf. \S\ref{pre:extended}).
When a splitting is chosen, or equivalently, $\TTC M$ is identified with $\TT M$
with an $H$-twisted Courant algebroid structure, we will use the notion 
\emph{$H$-twisted generalized} ($+$ \emph{structures}). Now, when the action of
$G$ preserves a splitting of $\TTC M$, then the reduced extended tangent bundle
in the theorem naturally splits and the twisting form on $Q$ can be explicitly
written down (cf. corollary \ref{courant:split}).

In investigating $T$-duality, we are guided by the detailed computation in
\cite{Hu} of the example of $\CC^2\setminus \{(0,0)\}$ with non-trivial twisting
class,
which we recall in \S\ref{torus:example}.
The following definition is crucial:

\vspace{0.1in}
\noindent
{\bf Definition \ref{double:biactdef}.} 
{\it 
Let $(\tmf g, \mf g, \hmf g)$ be the Manin triple defined by a Poisson Lie group
$G$ (cf. theorem \ref{app:liedouble}, also \cite{Lu}), with dual group $\hat G$.
The (infinitesimal) action of $\tmf g$ on $M$ is \emph{bi-Hamiltonian} if it is
induced by a $\JJC_1$-Hamiltonian action of $G$ together with a
$\JJC_2$-Hamiltonian action of $\hat G$.
}

\vspace{0.1in}
Suppose that the Manin triple $(\tmf g, \mf g, \hmf g)$ is the Lie algebras of
the (local) double Lie group $(\tilde G, G, \hat G)$ (cf. theorem
\ref{app:liedouble}). 
We impose two sets of assumptions, on the group $G$ (assumption
\ref{double:assum} $(0)$) and on the action of $\tilde G$ (the rest of
assumption \ref{double:assum}). Our first result in this direction is the
\emph{factorizable reduction}:

\vspace{0.1in}
\noindent
{\bf Theorem \ref{double:courant}.} {\it Under assumption \ref{double:assum} and
suppose that the action of $\tilde G$ on $M_0$ is proper and free, then the
reduced space $\tilde Q = M_0 / \tilde G$ of a bi-Hamiltonian action of
factorizable Poisson Lie group admits a natural effective Courant algebroid
(definition \ref{courant:extcourant}).}

\vspace{0.1in}
\noindent
With further restrictions, i.e. the reduction exists with respect to either of
the actions of $G$ and $\hat G$ as given in theorem \ref{kahler:reduct}, the
factorizable reduction as in theorem \ref{double:courant} can be factored in two
ways, $M_0 \xto{/G} Q \xto{/\hat G} \tilde Q$ or $M_0 \xto{/\hat G} \hat Q
\xto{/G} \tilde Q$. We then propose

\vspace{0.1in}
\noindent
{\bf Definition \ref{torus:courant}.} {\it The extended K\"ahler structures on
$Q$ and $\hat Q$ are \emph{Courant dual} to each other.
}

\vspace{0.1in}
\noindent
We note that any of the groups $\tilde G, G$ or $\hat G$ could be non-abelian.
Thus we have a candidate for the \emph{non-abelian} duality \emph{with
background twistings}. 
The more stringent but natural assumption that $G$ and $\hat G$ commute in
$\tilde G$ implies that $\tilde G$ is in fact a torus $\tilde T$.
The choice of terminology in the above is supported by the following theorem
when the action of $\tilde T$ preserves a splitting of $\TTC M$:

\vspace{0.1in}
\noindent
{\bf Theorem \ref{torus:duality}.} {\it After applying a natural
$B$-transformation on $M$, which does not change the reduced Courant algebroid
on $\tilde Q$, the twisting forms $h$ and $\hat h$ of the structures on $Q$ and
$\hat Q$ respectively satisfy:
$$\hat\pi^*\hat h - \pi^*h = d(\hat\Theta \wedge \Theta),$$
where $\pi$ and $\hat\pi$ are the quotient maps and $\Theta$ and $\hat \Theta$
are connection forms of principle torus bundles.
}

\vspace{0.1in}
\noindent
We point out that the equation above appears as part of the \emph{definition} of
$T$-duality with $H$-flux of principle torus bundles in the literature (also see
below). Here, it appears as a geometrical consequence. The notion of $T$-duality
group in the literature can be recovered (\S\ref{dual}) with our construction. 

We describe the content of the article in the following.
It's helpful to recall the basics of Lu's construction (see also
\S\ref{app:plie} appendix B). A Poisson Lie group $G$ is a Lie group with a
multiplicative Poisson structure, i.e. $m: G \times G \to G$ is a Poisson map.
Let $(M, \omega)$ be a symplectic manifold, the action of $G$ on $M$ is called
Poisson if the map $G\times M \to M$ defining the action is Poisson, with the
product Poisson structure on $G \times M$. In \cite{Lu}, Lu defined momentum
mapping for such Poisson actions (see also definition \ref{app:poissonmom},
theorem \ref{app:momequi}) and went on to show that symplectic reduction can be
carried out for Poisson actions with momentum mapping, although in general, the
symplectic structure $\omega$ is not invariant under Poisson actions.

The section \S\ref{recall} recalls the useful facts concerning the action of the
group of generalized symmetries $\tilde {\ms G} = \Diff(M) \ltimes \Omega^2(M)$,
the $H$-twisted Lie bracket on $\ms X = \Gamma(TM) \dsum \Omega^2_0(M)$, Courant
algebroid and generalized complex structures and explain in more detail the
notion of extended structures. These results are not new and details may be
found in, for example \cite{Cavalcanti, Gualtieri, Hitchin, Hu} and the
references therein.

We show, in \S\ref{poisson}, 
that the momentum mapping as defined in \cite{Lu} can be extended to the
generalized geometry (definition \ref{poisson:ham}), and the reduction
construction for symplectic manifold can be extended to generalized complex
manifold (theorem \ref{poisson:thm}), as well as generalized K\"ahler manifold
(theorem \ref{kahler:reduct}). 
Along the way, we obtain lemma \ref{subg:inv}, which can be viewed as an 
extension of Moser's argument for symplectic geometry (remark \ref{poisson:inv}).
We note that similar to the case of symplectic
geometry in \cite{Lu}, the generalized complex structure may not be preserved by
the group action. In fact, in our construction of generalized K\"ahler
reduction, none of the two generalized complex structures need to be preserved
by the group action, as long as certain subbundle of $\TT M$ is preserved
(remark \ref{kahler:nonpreserve}). We remark that reduction of Courant algebroid
(\S\ref{courant} appendix A) as well as reduction of generalized K\"ahler
structure have been discussed in various other works \cite{Gualtieri, Lin, Lin2,
Stienon}.

One of the features of generalized K\"ahler geometry is that the two generalized
complex structures are on the same footing, which is not at all obvious in the
classical K\"ahler geometry. In fact, this is one of the reasons that
generalized K\"ahler geometry could serve as the natural category of discussing
duality. Generalized K\"ahler geometry is relevant also from the result in
\cite{Gualtieri}, that it is equivalent to the bi-hermitian geometry, which has
been shown to be the string background for $N = (2,2)$ supersymmetry
(\cite{Gates}, \cite{Bredthauer} and references therein).
The notion of $T$-duality with $H$-flux in abelian case is proposed in
\cite{Bouwknegt} and then has been worked to much more general situations which
involve non-commutative \cite{Mathai} and non-associative \cite{Bouwknegt1}
geometries. The motivation in physics is that the physical theories on $T$-dual
spaces are isomorphic and thus provides insights to what the physics is about.
Here we concentrate on the more geometrical duality and leave the non-classical
cases to furture work.

We first describe the construction of $T$-duality with $H$-flux from the
existing literature in the following. To simplify matters, we restrict to
$T=S^1$, where many complications do not arise.
Let $p: E \to M$ be an $S^1$-principle bundle with connection form $\Theta \in
\Omega^1(E)$ and curvature form $\Omega \in \Omega^2(M)$. Let $H \in
\Omega^3(E)^{S^1}$ be a closed $S^1$-invariant $3$-form representing integral
class $[H] \in H^3(E, \ZZ)$. By construction, there is a form $h \in
\Omega^3(M)$ so that $p^*h = H - \Theta \wedge \hat\Omega$. Let $\hat\Omega \in
\Omega^2(M)$ be the integration of $H$ along the fiber of $E$, then
$[\hat\Omega] \in H^2(M, \ZZ)$ and there is a principle $S^1$-bundle $\hat p :
\hat E \to M$ whose first Chern class is $[\hat \Omega]$. In particular, we may
choose a connection form $\hat \Theta \in \Omega^1(\hat E)$ whose curvature form
is $\hat \Omega$. Let $\hat H = \hat p^*h + \hat \Theta \wedge \Omega$, then
$\hat H \in \Omega^3(\hat E)^{S^1}$ is closed and the pair $(\hat E, \hat H)$ is
said to be $T$-dual to the pair $(E, H)$. One may also consider the
correspondence space $E\times_M \hat E$, whose projection to $E$ and $\hat E$ is
denoted $\pi$ and $\hat \pi$ respectively. Then the forms satisfy $\hat \pi^*
\hat H - \pi^*H = d(\hat \Theta \wedge \Theta)$. We may summarize this
description with the following diagram:
$$\text{
\xymatrix{
& {E\times_M \hat E}\ar[dl]_{\pi} \ar[dr]^{\hat \pi} & \\
{(E, H; \Theta)} \ar[dr]_{p} & & {(\hat E, \hat H; \hat \Theta)} \ar[dl]^{\hat
p}\\
& {(M; h, \Omega, \hat \Omega)} &
}
}$$
For higher dimensional torus, it's argued (see \cite{Mathai}, \cite{Bouwknegt1}
and references therein) that various conditions are needed, on the action and
twisting form $H$, in order for the dual space to be classical. Otherwise, it
would be of one of the non-classical geometries.

The idea of applying generalized geometry in describing $T$-duality is
introduced by Gualtieri \cite{Gualtieri} and Cavalcanti \cite{Cavalcanti}, where
the first efforts were made. The guiding example for us is the example in
\S\ref{torus:dualtrans}. By this example, we see that it's possible for the
same function to serve as moment map for Hamiltonian group actions with respect
to either generalized complex structure and thus provides a diagram similar to
the one above. Another important input is from \cite{Cavalcanti}, where
Cavalcanti showed that the Courant algebroids defined by invariant sections on
$T$-dual $S^1$-principle bundles are isomorphic.

On the physics side, there is vast literature on $T$-duality, both with or
without $H$-flux, abelian or non-abelian, for principle bundles or fibration
with singular fibers. The approach of realizing dual theories by quotient
construction appeared in \cite{Rocek, Kiritsis}, where it's argued that gauging
different chiral currents produces dual $\sigma$-models. More recently, there is
work of Hull \cite{Hull}, which discusses $T$-duality in the \emph{doubled
formalism}. The formalism is to look at the correspondence space as principle
bundle of a doubled torus, consisting of the product of a dual pair of torus
with the natural pairing on the Lie algebra. Then the group automorphisms
preserving the pairing corresponds to the $T$-duality group. The idea of looking
to Poisson Lie group in considering duality goes back to a series of papers by
Klim\v c\'ik and/or \v Severa starting with \cite{Klimcik1, Klimcik2, Klimcik3},
where Poisson Lie target space duality was proposed as the framework of
non-abelian $T$-duality. The papers \cite{Klimcik4, Hull, Parkhomenko, Unge} and
references therein contain more recent developement in this direction.

Starting from \S\ref{double}, we discuss $T$-duality with $H$-fluxes in the
context of generalized (K\"ahler) geometry, which includes both abelian and
non-abelian groups.
In \S \ref{double}, we define the notion of bi-Hamiltonian action (definition
\ref{double:biactdef}) and discuss reduction of bi-Hamiltonian action of
factorizable Poisson Lie groups (theorem \ref{double:courant}). The main point
is that the reduced structure is an effective Courant algebroid on the reduced
space (definition \ref{courant:extcourant}). We note that the reduction
considered in \S\ref{double} can be factorized in two ways and in \S\ref{torus}
we define the two intermediate stages as being \emph{Courant dual} to each other
(definition \ref{double:basic}). Our construction then provides an isomorpism of
Courant algebroids defined by the invariant sections of Courant dual structures
(proposition \ref{torus:courantiso}), extending the result in \cite{Cavalcanti},
with a more geometrical method. The upshot is that in theorem
\ref{torus:duality}, we show that $T$-duality, as described above, can arise
from a special case of Courant duality. The notion of $T$-duality group is
essential in the full picture of $T$-duality with $H$-fluxes and we discuss it
in \S\ref{dual}. We note that it's more desirable that $T$-duality is
constructed starting from $(E, H; \Theta)$ instead of from the correspondence
space as the approach here. The construction of the correspondence space from
one of the reduced space will be discussed in \cite{Hu2}. 

To make the paper more self-contained, in \S\ref{courant} appendix A, we present
a construction of reduction of extended tangent bundles which is used in this
article. In \S\ref{app:plie} appendix B, we collect various facts on Lie
bialgebra, Poisson Lie group and Hamiltonian action.

\vspace{0.1in}
\noindent
{\bf Acknowledgement.} I would like to thank Francois Lalonde and DMS in
Universit\'e de Montr\'eal for support and excellent working conditions. I'd
like to thank Bernardo Uribe for stimulating discussions and correspondences. I
thank Vestislav Apostolov, Octav Cornea, Alexander Ivrii, Francois Lalonde, Sam
Lisi, and Alexander Cardona for helpful discussions. Special thanks go to my
family, for their support and understanding.

\section{Preliminaries}\label{recall}
We recall the preliminaries of generalized geometry and symmetries. As mentioned 
in the introduction, the results are not new and for details, the readers are 
referred to the literatures, for example \cite{Cavalcanti, Gualtieri, Hitchin, Hu},
and the references therein.
\subsection{} For a smooth manifold $M$, let $\TT M = TM \dsum T^*M$ and
$\tilde{\ms G} = \Diff(M) \ltimes \Omega^2(M)$. Let $\lambda, \mu \in \Diff(M)$
and $\alpha, \beta \in \Omega^2(M)$, then the product on $\tilde{\ms G}$ is
given by
$$(\lambda, \alpha)\cdot(\mu, \beta) = (\lambda\mu, \mu^*\alpha + \beta).$$
Let $\mf X = X + \xi$ with $X \in TM$ and $\xi \in T^*M$, then the (left) action
of $\tilde{\ms G}$ on $\TT M$ is given by
$$(\lambda, \alpha) \circ (X + \xi) = \lambda_* X + (\lambda^{-1})^*(\xi +
\iota_X \alpha).$$
The Lie algebra of $\tilde {\ms G}$ is $\tilde {\ms X} = \Gamma(TM) \dsum
\Omega^2(M)$ with the following Lie bracket:
$$[(X, A), (Y, B)] = ([X, Y], \LLC_X B - \LLC_Y A).$$
The $1$-parameter subgroup generated by $(X, A)$ is given by
$$e^{t(X, A)} = (\lambda_t, \alpha_t) = \left(e^{tX}, \int_0^t \lambda_s^* A
ds\right).$$
Following the above notation, for $B \in \Omega^2(M)$, we use $e^B$ to denote
the so-called $B$-transformation 
$$e^B\circ (X + \xi) = X + \xi + \iota_X B.$$

\subsection{} Let $H \in \Omega^3_0(M)$, i.e. $d H = 0$. The $H$-twisted Loday
bracket on $\TT M$ is defined by
$$(X + \xi) *_H (Y + \eta) = [X, Y] + \LLC_X \eta - \iota_Y (d\xi - \iota_X
H).$$
Let $\<X+ \xi, Y + \eta\> = \frac{1}{2}(\iota_X \eta + \iota_Y \xi)$, then $(\TT
M, *_H, \<,\>, a)$ defines a structure of Courant algebroid, with $a : \TT M \to
TM$ the natural projection (cf. definition \ref{app:courant} below).
The Loday bracket is not skew-symmetric, indeed we have
$$(X + \xi) *_H (Y + \eta) + (Y+\eta) *_H (X+\xi) = d\<X+\xi, Y + \eta\>.$$
The subgroup $\ms G = \Diff(M) \ltimes \Omega^2_0(M)$ is the group of symmetries
of the Courant algebroid structure with $H = 0$. The Lie algebra of $\ms G$ is
$\ms X = \Gamma(TM) \dsum \Omega_0^2(M)$ with the induced bracket. Let
$\tilde{\ms G}_H \subset \tilde{\ms G}$ be the symmetry group of the Courant
algebroid structure for general $H$ and $\tilde{\ms X}_H$ be its Lie algebra.
Consider the linear isomorphism:
$$\psi_H : \tilde{\ms X} \to \tilde{\ms X} : (X, A) \mapsto (X, A + \iota_X
H),$$
and the $H$-twisted Lie bracket
$$[,]_H : \ms X \times \ms X \to \ms X : [(X, A), (Y, B)]_H = ([X, Y], \LLC_X B
- \LLC_Y A + d\iota_Y \iota_X H),$$
then we have
\begin{prop}\label{recall:transl}\cite{Hu} For $H, H' \in \Omega^3_0(M)$,
$$[\psi_H(X, A), \psi_H(Y, B)]_{H+H'} = \psi_H[(X, A), (Y, B)]_{H'}$$
and $\psi_H : (\tilde\XXS_H, [,]) \to (\XXS, [,]_H)$ is Lie algebra isomorphism.
\qed
\end{prop}
\noindent
Let $\mf X = X + \xi \in \Gamma(\TT M)$, then $(X, d\xi) \in \ms X$ and
generates a $1$-parameter subgroup in $\tilde{\ms G}_H$:
$$e^{\psi_H^{-1}(X, d\xi)} = (\lambda_t, \alpha_t) = \left(e^{tX}, \int_0^t
\lambda_s^* (d\xi - \iota_X H) ds\right).$$
The infinitesimal action of $\mf X$ on $\mf Y \in \Gamma(\TT M)$ that generates
the above subgroup is:
$$\mf X \circ_H \mf Y = - \mf X*_H \mf Y.$$

\subsection{} Let $\JJ : \TT M \to \TT M$ be a generalized almost complex
structure on $M$, that is, $\JJ^2 = -\1$ and $\JJ$ is orthogonal with respect to
the pairing $\<,\>$. Let $L \subset \TT_\CC M$ be the $i$-eigensubundle of
$\JJ$, then $L$ is isotropic and $\JJ$ defines an $H$-twisted generalized
complex structure if $L$ is involutive with respect to the $H$-twisted Loday
bracket $*_H$.
Examples of generalized complex structures include the symplectic and complex
structures. Let $\omega$ (resp. $J$) be a symplectic (resp. complex) structure
on $M$, then the corresponding generalized complex structure is defined by the
respective isotropic subbundles:
$$L_\omega = \{X - i\iota_X\omega | X \in TM\} \text{ and } L_J = \{X+\xi + i
(J(X) - J^*(\xi)) | X \in TM, \xi \in T^*M\}.$$

\subsection{}\label{spinors} The space of complex valued differential forms
$\Omega^\bullet(M; \CC)$ is the spinor space of generalized geometry. Let $d_H =
d - H\wedge$ be the $H$-twisted differential on $\Omega^\bullet(M; \CC)$. Each
maximally isotropic subbundle $L \subset \TT_\CC M$ corresponds to a pure line
subbundle $U$ of $\wedge^\bullet T^*_\CC M$:
$$U = \Ann_C(L) := \{\rho \in \wedge^\bullet T^*_\CC M| \mf X \cdot \rho =
\iota_X \rho + \xi \wedge \rho = 0 \text{ for all } \mf X = X + \xi \in L\},$$
where $\cdot$ stands for the Clifford multiplication.
A (nowhere vanishing) local section $\rho$ of $U$ is called a \emph{pure spinor}
associated to the subbundle $L$. The integrability of $L$ with respect to the
$H$-twisted Courant bracket is equivalent to the condition $d_H(\Gamma(U))
\subset \Gamma(U_1)$, where $U_1 = \Gamma(\TT_\CC M) \cdot U$ via Clifford
multiplication. More explicitly, there is a unique local section $\mf Y = Y +
\eta$ of $\bar L$, so that 
\begin{equation}\label{subg:integ}
d_H \rho = d\rho - H \wedge \rho = \mf Y \cdot \rho = \iota_Y \rho + \eta \wedge
\rho,
\end{equation}
where we use the convention of $d_H$ as in \cite{Mathai1}. For a generalized
complex structure $\JJ$, the complex line bundle $U$ is called the
\emph{canonical bundle} of $\JJ$.

\subsection{} Via \cite{Kosmann3}, we have the following definition of a Courant
algebroid:

\begin{defn}\label{app:courant}
Let $E \to M$ be a vector bundle. A \emph{Loday bracket} $*$ on $\Gamma(E)$ is a
$\RR$-bilinear map satisfying the Jacobi identity, i.e. for all $\mf X, \mf Y,
\mf Z \in \Gamma(E)$, \begin{equation}\label{app:courantdefn1}
\mf X * (\mf Y * \mf Z) = (\mf X*\mf Y) * \mf Z + \mf Y*(\mf X * \mf Z).
\end{equation}
$E$ is a \emph{Courant algebroid} if it has a \emph{Loday bracket} $*$ and a
non-degenerate symmetric pairing $\<,\>$ on the sections, with 
an \emph{anchor map} $a : E \to TM$ which is a vector bundle homomorphism so
that
\begin{eqnarray}
\label{app:courantdefn2} a(\mf X) \<\mf Y, \mf Z\> & = & \<\mf X, \mf Y*\mf Z +
\mf Z* \mf Y\> \\
\label{app:courantdefn3} a(\mf X) \<\mf Y, \mf Z\> & = & \<\mf X* \mf Y, \mf Z\>
+ \<\mf Y, \mf X* \mf Z\>.
\end{eqnarray}
\end{defn}
The skew-symmetrization $[,]$ of $*$ in the definition is also called the
Courant bracket.
In particular, the datum $(\TT M, *_H, \<,\>, a)$ as given in previous
subsections, for $H \in \Omega^3_0(M)$, are examples of Courant algebroid, where
the corresponding Courant bracket is usually denoted $[,]_H$.

\subsection{}\label{pre:extended} Let $\TTC M$ be a Courant algebroid which fits
into the following extension:
$$0 \to T^*M \to \TTC M \xto{a} TM \to 0,$$
so that $a$ is the anchor map. Such Courant algebroid is called \emph{exact}
\cite{Severa}. The set of isotropic splitting $s : TM \to \TTC M$ is non-empty
and is a torsor over $\Omega^2(M)$. The choice of such $s$ determines a form $H
\in \Omega^3_0(M)$ and $\TTC M$ can then be identified with the datum $(\TT M,
*_H, \<,\>, a)$ as discussed above. The action of $B \in \Omega^2(M)$ on the set
of splittings translates into $H \mapsto H + dB$ on the corresponding forms. It
follows that $[H] \in H^3(M; \RR)$ is well-defined and is the \emph{\v Severa
class} of $\TTC M$. We use the notion \emph{extended} ($+$ \emph{structures}) to
emphasize the absence of a splitting while reserve \emph{twisted generalized}
for the situation where a splitting is (or can be explicitly) chosen. For
example, an extended complex structure $\JJC$ will represent a twisted
generalized complex structure $\JJ$ on $\TT M$ (once a splitting is chosen),
which is integrable with respect to a twisted Loday bracket $*_H$, where $[H]$
gives the \v Severa class of the extended tangent bundle $\TTC M$ defined by the
Courant algebroid structure $(\TT M, *_H, \<,\>, a)$. Given a different choice
of splitting of $\TTC M$, $\JJC$ will represent $\JJ_B$, which is $\JJ$
transformed by some $B \in \Omega^2(M)$ and is integrable with respect to
$*_{H+dB}$ on $\TT M$. We note that the Courant algebroids are \emph{identical}
(not only isomorphic) in either cases, since the difference is only the choice
of a splitting that gives the identification to $\TT M$.

\section{Poisson Lie actions and reductions} \label{poisson}
In this section, we describe the reduction of generalized complex and K\"ahler
manifolds via the Hamiltonian action of Poisson Lie groups. This extends the
reduction construction of \cite{Hu} for Hamiltonian action of Lie groups and
that of \cite{Lu} for Poisson Lie action on symplectic manifolds, which we
describe in the appendix B (\S\ref{app:plie}). Again, when we use $\TT M$, $\JJ$ and
etc, we assume a choice of splitting of the extended tangent bundle $\TTC M$ and
identify the corresponding structures as $H$-twisted generalized structures.

\subsection{} 
We first discuss the invariance of $\JJ$ under generalized actions. Direct
computation shows
\begin{lemma}Let $(\lambda, \alpha) \in \tilde\GGS$ and $\rho$ be the pure
spinor defining $\JJ$, then $(\lambda, \alpha) \circ \rho :=
(\lambda^{-1})^*(e^{-\alpha}\rho)$ is the pure spinor defining $(\lambda,
\alpha) \circ \JJ$. If $\JJ$ is $H$-twisted integrable, then $(\lambda,
\alpha)\circ \JJ$ is $(\lambda,\alpha)\circ H $-twisted integrable, where
$(\lambda,\alpha)\circ H = (\lambda^{-1})^*(H - d\alpha)$. We have $d_{(\lambda,
\alpha)\circ H} (\lambda, \alpha)\circ \rho = (\lambda, \alpha)\circ d_H\rho$.
\qed
\end{lemma}
\begin{remark}\label{subg:cplx}
\rm{We note that when considering generalized symmetries, we do not have to
restrict to real $2$-forms to stay with \emph{real} twisting form, e.g. the
group $\Diff(M) \ltimes (\Omega^2(M) \dsum i \Omega^2_0(M))$ acts on $\TT_\CC
M$. The infinitesimal action of $(X, A) \in \tilde\XXS \dsum i \Omega^2_0(M)$ on
the spinors is then given by
$$(X, A) \circ \rho = -\LLC_X \rho - A\wedge \rho.$$
For $\mf X = X + \xi \in \Gamma(\TT_\CC M)$ so that $X \in \Gamma(TM)$, let $(X,
A) = (X, d\xi - \iota_X H)$ and we compute the infinitesimal action on a section
$\rho$ of the canonical bundle of $\JJ$:
\begin{equation}\label{gency:commu}
\mf X \circ_H \rho = (-d_H + \mf Y \cdot) \mf X \cdot \rho - \<\mf X, \mf Y\>
\rho.
\end{equation}
We caution that when a generalized complex structure is concerned, such
\emph{complex} actions in general might not preserve the real index.
}\end{remark}
\begin{lemma}\label{subg:inv}
Suppose that $L$ defines the extended complex structure $\JJC$ and $\mf X_t =
X_t + \xi_t \in \Gamma(L \inter \TTC M \dsum iT^* M)$ is a family of sections parametrized by $\RR$.
Let $\tilde \lambda_t = (\lambda_t, \alpha_t)$ be the family of generalized symmetries generated by $\mf X_t$.
Suppose that for each $p \in M$ there is an open neighbourhood $U_p$ and a
compact set $V_p$ so that $\{\lambda_t \circ U_p\} \subset V_p$ for all $t$. Then
$\tilde \lambda_t$ preserves $\JJC$ for all $t$.
\end{lemma}
{\it Proof:} Choose a splitting and identify the structures with $H$-twisted
structures. Starting from any $p \in M$ and $t_0 \in \RR$. Suppose that 
$\rho_{t_0} = \rho$ is a local section of the canonical bundle $U$
of $\JJ$ and $\rho_t = (\lambda_t, \alpha_t)^* \rho := (\lambda_t,
\alpha_t)^{-1} \circ \rho = e^{\alpha_t}\lambda_t^* \rho$. Then $\rho_t$ is 
a local section of the canonical bundle $U_t$ of $\JJ_t = (\lambda_t, 
\alpha_t)^{-1} \circ \JJ$. Direct computation shows that
$$\frac{d}{dt} \rho_t = \left.\frac{d}{ds}\right|_{s = 0} \rho_{t+s} =
(\lambda_t, \alpha_t)^*((d_H - \mf Y \cdot) \mf X_t \cdot \rho + \<\mf X_t, \mf
Y\> \rho).$$
Then by assumption $\frac{d}{dt} \rho_t = f_t \rho_t$ for $f_t =
\lambda_t^*\<\mf X_t, \mf Y\>$. With the initial condition of $\rho_{t_0} =
\rho$ we get
$$\rho_t = e^{\int_{t_0}^t f_s ds} \rho.$$
It follows that $U_t = U$
wherever both $\rho$ and $\rho_t$ are defined, e.g. for a neighbourhood of $p$.
Since $t_0$ is arbitrary and by compactness 
assumption, we see that $U_t = U$ for all $t$.

The argument above shows that $L$ is preserved by the family of symmetries
generated by $\mf X_t$, which is independent of the splitting chosen. The
proposition then follows.
\qed

\begin{remark}\label{poisson:inv}
\rm{ 
Of course, when $M$ is compact, the condition in the proposition automatically
holds. 
From the proof, we also see that when $d_H \rho = 0$, not only the canonical
line bundle is preserved, the spinor $\rho$ is preserved as well.
}
\end{remark}
\begin{remark}\label{poisson:moser}
\rm{
The Moser's argument in symplectic geometry can be seen as a special case of the
above lemma. The Moser's argument goes as following (see \cite{McDuff}).
Consider a smooth family of symplectic forms $\omega_t = \omega_0 + d\beta_t$
and $\eta_t = \frac{d}{dt} \beta_t$. Let $Y_t$ be defined by $\iota_{Y_t}
\omega_t + \eta_t = 0$ and $\phi_t$ be the family of diffeomorphisms generated
by $Y_t$ via $\frac{d}{dt} \phi_t = \phi_{t*} (Y_t)$, then $\phi_t^* \omega_t =
\omega_0$.

In light of lemma \ref{subg:inv}, we consider $\varphi_t =
\phi_t^{-1}$, which is generated by the family of vector fields $X_t = -
\varphi_{t*} (Y_t)$. Then we define $\xi_t = \iota_{X_t} \omega_0 =
-\phi_t^*(\eta_t)$. It follows that $\mf X_t = X_t - i \xi_t \in
\Gamma(L_{\omega_0})$. The lemma then implies that the following family of
symmetries preserves $L_{\omega_0}$:
$$(\varphi_t, \alpha_t) = \left(\varphi_t, -id\int_0^t \varphi_s^*\xi_s
ds\right) = \left(\varphi_t, -id \int_0^t \eta_s ds\right).$$
In this case, we have $\rho_0 = e^{i\omega_0}$ and $d\rho_0 = 0$. The proof of
the lemma then implies that $\rho_0$ is preserved:
$$e^{i\omega_0} = (\varphi_t, \alpha_t)^*e^{i\omega_0} = e^{-id \int_0^t \eta_s
ds} \varphi_t^* (e^{i\omega_0}) = e^{-id\beta_t} (\phi_t^*)^{-1}
(e^{i\omega_0}),$$
which is equivalent to $\phi_t^* \omega_t = \omega_0$ as in Moser's argument.
}
\end{remark}

\subsection{}
We will use the following conventions:
\begin{convention}\label{app:conv}
Given a Lie group $G$, the Lie algebra $\mf g$ of $G$ is identified as the
tangent space $T_e G$ at identity, as well as the space of right invariant
vector fields, i.e. $\tau \to X_\tau(g) = (R_g)_* \tau$. Then the dual $\hmf g 
= \mf g^*$ of the Lie algebra is identified with the space of right invariant 
$1$-forms on $G$. Let $\theta^r_{\hat \tau} \in \Omega^1(G)$ denote the right invariant 
$1$-form on $G$ with $\theta^r_{\hat \tau}(e) = \hat \tau$ and $\theta^l_{\hat \tau}$ the left invariant
$1$-form on $G$ with $\theta^l_{\hat \tau}(e) = \hat \tau$, for $\hat \tau \in \hmf g$. Given a
Poisson manifold $P$ with Poisson tensor $\pi_P$, we consider $\pi_P$ also as a
map $\pi_P : TP^* \to TP$ defined by $\iota_{\pi_P(\xi)}\eta = \pi_P(\xi, \eta)$
for $\xi, \eta \in \Omega^1(P)$.
\end{convention}
 We note that for $\tau \in \mf g$, the right invariant vector field $X_\tau$
generates left action on $G$ by the $1$-parameter subgroup $g_t = e^{t\tau}$.
Thus the left action of $G$ on $M$ induces homomorphism of Lie algebras $\xi
\mapsto X^M_\tau$ where $X^M_\tau$ is the infinitesimal action generated by $\tau$,
while the right action of $G$ induces \emph{anti}-homomorphism of Lie algebras.
With this convention, the map $\pi_P$ and the Lie algebra (anti)-homomorphism
are opposite to the convention used in \cite{Lu} and \cite{Lu1}. 
In the following, we will only consider \emph{left} actions. 
We collected the relevant
definitions and results on Poisson Lie groups and actions in the appendix
(\S\ref{app:plie}).
\subsection{} \label{poisson:induce}
The basic setup is the following. Let $(M, \JJC)$ be an extended complex
manifold with extended tangent bundle $\TTC M$ and anchor $a: \TTC M \to TM$.
Then there is a natural induced Poisson structure $\pi_\JJC$ on $M$ defined by
$$\pi_\JJC : T^*M \to \TTC M \xto{\JJC} \TTC M \xto{a} TM.$$
Let $(G, \pi_G)$ be a Poisson Lie group with Poisson structure $\pi_G$. Let
$\sigma : G \times M \to M$ be a (left) Poisson action with equivariant moment
map $\mu : M \to \hat G$ (cf. definition \ref{app:poissonmom}).
If a splitting is chosen, we use $\TT M$, $\JJ$ and $H$-twisted when referring
to the respective structures. We show
\begin{lemma}\label{poisson:moment}
Suppose that $G$ is connected. 
Let $\mu : M \to \hat G$ be an equivariant moment map as in definition
\ref{app:poissonmom} and $\JJ(\mu^*\hat\theta_\tau) = \mf X_\tau = X_\tau +
\xi_\tau $ for $\tau \in \mf g$, then
\begin{equation*}
\iota_{X_\tau} \mu^*(\hat\theta_\tau) = \iota_{X_\tau} \xi_\tau = 0 \text{ and }
[(X_\tau, d\xi_\tau), (X_\omega, d\xi_\omega)]_H = (X_{[\tau, \omega]},
d\xi_{[\tau, \omega]}).
\end{equation*}
\end{lemma}
{\it Proof: }
First $\iota_{X_\tau}\mu^*\hat \theta_\tau = \mu^*(\iota_{\hat\XXC_\tau}
\hat\theta_\tau) = \mu^*(\hat\pi_G(\hat \XXC_\tau, \hat \XXC_\tau)) = 0$, where
$\hat \XXC_\tau$ is the dressing vector field generated by $\tau \in \mf g$ (cf.
definition \ref{app:dressing}). We then compute 
\begin{equation*}\begin{split}
& [\JJ(\mu^*\hat\theta_\tau) + i\mu^*\hat\theta_\tau,
\JJ(\mu^*\hat\theta_\omega) + i\mu^*\hat\theta_\omega]_H \\
= & [X_\tau, X_\omega] + \LLC_{X_\tau}\xi_\omega - \iota_{X_\omega}d\xi_\tau +
\iota_{X_\omega}\iota_{X_\tau} H + i(\LLC_{X_\tau}\mu^*\hat\theta_\omega -
\iota_{X_\omega}d\mu^*\hat\theta_\tau),\\
\end{split}
\end{equation*}
and the imaginary part is
\begin{equation*}
\begin{split}
& \LLC_{X_\tau}\mu^*\hat\theta_\omega - \iota_{X_\omega}d\mu^*\hat\theta_\tau =
\mu^*(\LLC_{\mu_*X_\tau}\hat\theta_\omega -
\iota_{\mu_*X_\omega}d\hat\theta_\tau) \\
= & \mu^*(\LLC_{\hat \XXC_\tau}\hat\theta_\omega - \iota_{\hat
\XXC_\omega}d\hat\theta_\tau)
= \mu^*([\hat\theta_\tau, \hat\theta_\omega]^*) = \mu^*\hat\theta_{[\tau,
\omega]}
\end{split}
\end{equation*}
Thus $X_{[\tau, \omega]} = [X_\tau, X_\omega]$ and $\xi_{[\tau, \omega]} =
\LLC_{X_\tau}\xi_\omega - \iota_{X_\omega}d\xi_\tau +
\iota_{X_\omega}\iota_{X_\tau} H$. The lemma follows.
\qed

We note that from above lemma, the symmetry generated by $\mf X_{[\tau,
\omega]}$ coincides with that of $\mf X_\tau *_H \mf X_\omega$, which only
depends on the Loday bracket.
\begin{defn}\label{poisson:ham}
The action of a Poisson Lie group $G$ on an $H$-twisted generalized complex
manifold $(M, \JJC)$ is \emph{Hamiltonian} with \emph{moment map} $\mu : M \to
\hat G$, if the action is Poisson with respect to $\pi_\JJC$, 
$\mu$ is an equivariant moment map as in definition \ref{app:poissonmom}, so
that the $G$-action on $\TTC M$ is generated by $\JJC(\mu^*\hat\theta) = \mf
X_\mu$, via the Loday bracket $*$.
\end{defn}
\begin{remark}\label{poisson:pres}
\rm{ We recall that in the Poisson category, the Poisson action of a Poisson Lie
group does not have to preserve the Poisson structure. Thus the action as
defined above does not have to preserve the extended complex structure $\JJC$.
The lemma \ref{subg:inv} implies that the action on $\TTC M$ generated by
$\JJC(\mu^*\hat\theta) + i\mu^*\hat\theta$ does preserve the structure $\JJC$.
Thus the non-invariance under the action above can be seen as due to the
non-closedness of $\hat \theta$.
When the Poisson structure on $G$ is trivial, we have the definition for
Hamiltonian actions of Lie groups \cite{Hu}. By theorem \ref{app:momequi}, $\mu$
is a Poisson map. Let $M_0 = \mu^{-1}(\hat e)$, then $\mu_*(\pi_\JJC|_{M_0}) =
\pi_{\hat G}|_{\hat e} = 0$, and $X_\mu$, i.e. the geometrical action of $G$,
preserves $M_0$.
}
\end{remark}

\subsection{} \label{poisson:assum} We may consider reduction by Hamiltonian
Poisson Lie group action. Assume that 
\begin{enumerate}
\item The identity $\hat e\in \hat G$ is a regular value of $\mu$,
\item (\emph{the geometrical part of}) $G$ acts freely on $M_0$.
\end{enumerate}
\begin{lemma}\label{poisson:equiv}
The restriction of $(\mu^*\hat \theta)$, $\JJC(\mu^*\hat\theta)$ and $L \dsum
(\mu^*\hat \theta)$ are $G$-equivariant subbundles.
\end{lemma}
{\it Proof: }
Choose a splitting. It's enough to show that the infinitesimal actions preserve
the subbundles:
\begin{equation*}\begin{split}
(X_\omega, d\xi_\omega) *_H (\mu^*\hat\theta_\tau) & =
\LLC_{X_\omega}\mu^*\hat\theta_\tau = \mu^*(\LLC_{\hat  \XXC_\omega}
\hat\theta_\tau) = \mu^*\hat\theta_{[\omega, \tau]}, \\
(X_\omega, d\xi_\omega) *_H (X_\tau + \xi_\tau) & = [X_\omega, X_\tau] +
\LLC_{X_\omega}\xi_\tau - \iota_{X_\tau}d\xi_\omega +
\iota_{X_\tau}\iota_{X_\omega} H = X_{[\omega, \tau]} + \xi_{[\omega, \tau]}, \\
(X_\omega, d\xi_\omega) *_H (Y + \eta) & = [X_\omega, Y] + \LLC_{X_\omega}\eta -
\iota_Yd\xi_\omega + \iota_Y\iota_{X_\omega} H \\
& = (X_\omega + \xi_\omega + i\mu^*\hat \theta_\omega) *_H (Y+\eta) +
i\iota_Y\mu^*(d\hat\theta_\omega),
\end{split}
\end{equation*}
for $Y+\eta \in \Gamma(L)$. We note that $\iota_Y\mu^*(d\hat\theta_\omega) =
-\frac{1}{2}\iota_Y\mu^*([\hat\theta, \hat\theta]_\omega) \in (\mu^*\hat
\theta)$.
\qed
\begin{theorem}\label{poisson:thm}
Suppose that $G$ is compact and assumptions $(1)$ and $(2)$, there is a natural
extended complex structure on the quotient $Q = M_0 / G$. When the action of $G$
preserves a splitting of $\TTC M$, the reduced structure $\TTC Q$ admits a
natural splitting up to a choice of connection form of $M_0 \to Q$.
\end{theorem}
{\it Proof: }
By $\pi_{\hat G}|_{\hat e} = 0$ we compute on $M_0$:
$$\<\mu^*\hat\theta_\tau, \JJC(\mu^*\hat\theta_\omega)\> = \iota_{X_\omega}
\mu^*\hat\theta_\tau = \mu^*\iota_{\hat\XXC_\omega} \hat\theta_\tau =
\mu^*\pi_{\hat G}(\hat\XXC_\omega, \hat\XXC_\tau) = 0.$$

Then $\KKC = \mu^*\hat \theta, \KKC' = \JJC(\mu^*\hat\theta)$ satisfy the
conditions of lemma \ref{courant:reduce} $(1)$. Thus $\TTC_\mu M_0 =
\frac{\Ann(\mu^*\hat\theta, \JJC(\mu^*\hat\theta))}{(\mu^*\hat\theta,
\JJC(\mu^*\hat\theta))}$ descends to an extended tangent bundle $\TTC_\mu Q$ on
$Q$.
Consider $(L \dsum (\mu^*\hat \theta)) \inter \Ann(\mu^*\hat\theta,
\JJC(\mu^*\hat\theta))$, then it induces a subbundle $L_0$ in $\TTC_{\mu, \CC}
M_0$ which coincides with the image of $L$ under the subquotient. By lemma
\ref{poisson:equiv}, the bundle $L_0$ is $G$-equivariant and descends to a
subbundle $L_\mu$ of $\TTC_{\mu,\CC} Q$. That $L_\mu$ is maximally isotropic
with real index $0$ and integrable follows from the same properties of $L$. Thus
$L_0$ defines an extended complex structure $\JJC_\mu$. The last sentence
follows from corollary \ref{courant:split}.
\qed

\subsection{}
Let $(M, \JJC_1)$ be an extended complex manifold. A second extended complex
structure $\JJC_2$ makes $(M, \JJC_1, \JJC_2)$ into an \emph{extended K\"ahler}
manifold if they are both defined on the same extended bundle $\TTC M$ and $\GGC
= -\JJC_1\JJC_2 = -\JJC_2\JJC_1$ defines generalized metric (see
\cite{Hitchin3}) on $\TTC M$, i.e. $\<\GGC \cdot, \cdot\>$ defines a metric on
$\TTC M$. 
We show that just as symplectic reduction admits induced K\"ahler structure when
the original manifold is K\"ahler with $G$ preserving the complex structure,
generalized complex reduction with respect to $\JJC_1$ would admit extended
K\"ahler structure if  $\JJC_2$ is preserved.
\begin{defn}\label{kahler:ham} A Poisson action of Poisson Lie group $G$ on an
extended K\"ahler manifold $(M, \JJC_1, \JJC_2)$ is \emph{$\JJC_1$-Hamiltonian}
if it is Hamiltonian with respect to $\JJC_1$ and preserves $\JJC_2$. 
\end{defn}

\subsection{} We note that $\Ann(\mu^*\hat\theta, \JJC_1(\mu^*\hat\theta))$ is
not preserved by $\JJC_2$:
$$\Ann(\mu^*\hat\theta, \JJC_1(\mu^*\hat\theta)) \inter
\JJC_2(\Ann(\mu^*\hat\theta, \JJC_1(\mu^*\hat\theta))) = \Ann(\mu^*\hat\theta,
\JJC_1(\mu^*\hat\theta), \JJC_2(\mu^*\hat\theta), \GGC(\mu^*\hat\theta)).$$
The right hand side of above equation is again $G$-equivariant subbundle when
restricting to $M_0$, as the two terms on the left are both so. 
\begin{lemma}\label{kahler:linear}
Let $\VV = V\dsum V^*$ and $(\JJ_1, \JJ_2; \GG)$ be a linear generalized
K\"ahler structure. Given subspace $K \subset V^*$, let $\UU^j_K = \Ann(K,
\JJ_j(K))$, $\WW_K   = \Ann(K, \JJ_1(K), \JJ_2(K), \GG(K))$ and $(\UU^j_K)_\CC$,
$(\WW_K  )_\CC$ be the respective complexified versions, then
$L_j \inter (\WW_K  )_\CC = L_j \inter (\UU^l_K)_\CC \text{ for } j \neq l$.
If
\begin{enumerate}
\item $K + \JJ_1(K) \subset \UU^1_K$,
\end{enumerate}
then the following decomposition holds
$$\UU^1_K = \WW_K   \dsum (K + \JJ_1(K)).$$
The $+$ above becomes $\dsum$ if we suppose further that
\begin{enumerate}
\item[$(2)$] $\JJ_1(K) \inter V^* = \{0\}$.
\end{enumerate}
\end{lemma}
{\it Proof: }
Obviously $L_j \inter (\WW_K  )_\CC \subset L_j \inter (\UU^l_K)_\CC$. Let $\mf
X \in L_j \inter (\UU^l_K)_\CC$, then $\JJ_j(\mf X) = i \mf X$ and $\<\mf X, K\>
= \<\mf X, \JJ_l(K)\> = 0$. Then by orthogonality of $\JJ_j$ we find that $\<\mf
X, \JJ_j(K)\> = \<\mf X, \GG(K)\> = 0$, i.e. $\mf X \in L_j \inter (\WW_K 
)_\CC$.

For any subspace $W \subset \VV$, let $W^\perp = \Ann(\GG W)$, then $\VV = W
\dsum \Ann(\GG W)$. In particular
\begin{equation}\label{kahler:decomp}
\begin{split}
\VV & = \WW_K   \dsum (K + \JJ_1(K) + \JJ_2(K) + \GG(K)) = \tilde \WW_K   \dsum
(\JJ_1(K) + \JJ_2(K) + \GG(K)) \\
& = \UU^j_K \dsum (\JJ_l(K) + \GG(K)) \text{ for } j \neq l,
\end{split}
\end{equation}
where $\tilde \WW_K   = \Ann(K, \JJ_1(K), \JJ_2(K))$.
With condition $(1)$, by the last expression in \eqref{kahler:decomp} for $j =
1$ and $l = 2$, we see that the decomposition in the statement holds. With the
condition $(2)$, we have $K + \JJ_1(K) = K \dsum \JJ_1(K)$ and it follows that
all $+$'s in \eqref{kahler:decomp} are $\dsum$'s.
\qed
\begin{lemma}\label{kahler:sequence}
Continue from lemma \ref{kahler:linear} and let $N = a\circ \JJ_1(K)$ where $a:
\VV \to V$ is the projection, then 
there is a self-dual exact sequence
$$0 \to W^*_K \xto{a^*_K} \WW_K   \xto{a_K} W_K \to 0 \text{ where } W_K =
\frac{\Ann_V(K)}{N} \text{ and } W_K^* = \frac{\Ann_{V^*}(N)}{K}.$$
The restriction $\<,\>_K$ of $\<,\>_K$ on $\WW_K  $ is non-degenerate pairing
and $(\JJ_1, \JJ_2; \GG)$ restricts to generalized K\"ahler  structure
$(\JJ_{1,K}, \JJ_{2,K}; \GG_K)$ on $\WW_K  $ with respect to the pairing
$\<,\>_K$. 
The inclusion $\WW_K \into \Ann(K, \JJ_1(K))$ induces natural isomorphism $\WW_K
\simeq \VV_K$ in \S\ref{courant}, and the extension sequences correspond.
\end{lemma}
{\it Proof: } Note that $\WW_K  $ is preserved by $\GG$ we see that for any $\mf
X \in \WW_K  $ such that $\<\mf X, \WW_K  \> = 0$, it must satisfy $\<\mf X,
\GG(\mf X)\> = 0$, i.e. $\mf X = 0$. It implies that the restriction $\<,\>_K$
is nondegenerate.

Let $a_K : \WW_K   \to W_K$ be the map induced from the projection $a$. The
kernel of $\UU^1_K \to W_K$ is $\Ann_{V^*}(N) \dsum \JJ_1(K)$. It follows that
the kernel of $a_K$ is
$$\ker a_K = (\Ann_{V^*}(N) \dsum \JJ_1(K)) \inter \WW_K  .$$
Note that $\UU_K^1 = \VV_k \dsum (K \dsum \JJ_1(K))$ and $K\dsum \JJ_1(K)\subset
\Ann_{V^*}(N) \dsum \JJ_1(K)$, we find that $\Ann_{V^*}(N) \dsum \JJ_1(K) = \ker
a_K \dsum (K \dsum \JJ_1(K))$, thus $\ker a_K \simeq \frac{\Ann_{V^*}(N)}{K}$.
Now $\ker a_K$ is maximally isotropic with respect to $\<,\>_K$ and the
self-duality follows.
The last sentence follows from direct checking.
\qed

Similar to the classical K\"ahler case, we have:
\begin{theorem}\label{kahler:reduct}
Let the $G$ action on an extended K\"ahler manifold $(M, \JJC_1, \JJC_2; \GGC)$
be $\JJC_1$-Hamiltonian with moment map $\mu : M \to \hat G$. Suppose that the
assumptions in \S\ref{poisson:assum} hold, then there is a natural extended
K\"ahler structure on the quotient $Q =  M_0/G$. When the $G$-action preserves a
splitting of $\TTC M$, the reduced structure splits, up to a choice of
connection form on $M_0 \to Q$.
\end{theorem}
{\it Proof: } All the bundles in the proof will be on spaces at $\mu = \hat e$,
either level set or reduced space.
Let $\TTC'_\mu  M_0 = \Ann(\mu^*\hat\theta, \JJC_1(\mu^*\hat\theta),
\JJC_2(\mu^*\hat\theta), \GGC(\mu^*\hat\theta))$ be the subbundle of $\TTC M
|_{M_0}$, then it is a $G$-equivariant subbundle. 
Let $\TTC_\mu M_0$ be defined as in theorem \ref{poisson:thm} for $\JJC_1$. From
lemma \ref{kahler:sequence}, the bundles $\TTC'_\mu  M_0$ and $\TTC_\mu M_0$ are
naturally isomorphic via the inclusion of $\TTC'_\mu M_0$ in
$\Ann(\mu^*\hat\theta, \JJC_1(\mu^*\hat\theta))$.
We define a new bracket $*_1$ on $\Gamma(\TTC'_\mu  M_0)^G$ by the following:
$$\mf X *_1 \mf Y = \pi_1(\mf X * \mf Y), \text{ for } \mf X, \mf Y \in
\Gamma(\TTC'_\mu  M_0)^G \subset \Gamma(\Ann(\mu^*\hat\theta,
\JJC_1(\mu^*\hat\theta)))^G,$$
where $\pi_1$ is the projection $\Ann(\mu^*\hat\theta, \JJC_1(\mu^*\hat\theta))
\to \TTC'_\mu  M_0$ defined by the first decomposition in lemma
\ref{kahler:linear}. By definition, $\Gamma(\TTC'_\mu  M_0)^G$ is closed under
$*_1$. By construction, $\pi_1$ coincides with the projection
$\Ann(\mu^*\hat\theta, \JJC_1(\mu^*\hat\theta)) \to \TTC_\mu  M_0$ under the
natural isomorphism $\TTC_\mu  M_0 \simeq \TTC'_\mu  M_0$. Thus as in the
theorem \ref{poisson:thm}, 
the structure $(\TTC'_\mu  M_0, \<,\>_\mu, *_1)$ descends to an extended tangent
bundle $\TTC Q$ on $Q$. 
\qed
\begin{remark}\label{kahler:nonpreserve}
\rm{ We notice from the proof that, in order to have extended K\"ahler
reduction, even the extended complex structure $\JJC_2$ doesn't have to be
preserved by the $G$-action either. The only thing that needs to be preserved is
the intersection $L_2 \inter \TTC'_\mu  M_0$. Here, unlike the case in theorem
\ref{poisson:thm}, where $L_1 \dsum (\mu^*\hat\theta)$ being equivariant
provides descending of $\JJC_1$, lemma \ref{kahler:linear} implies that such
flexibility doesn't apply to $\JJC_2$.
}
\end{remark}
\begin{remark}\label{kahler:other}
\rm{Generalized K\"ahler reduction have been constructed by several other works,
e.g. \cite{Bursztyn}, \cite{Lin} and \cite{Stienon}, with various generalities.
The construction we describe here, which fits our needs for discussing duality,
has not appeared in the stated form. In particular, we allow non-trivial
$B$-field action and we only require the action preserve one of the generalized
complex structures.
}
\end{remark}

\section{Bi-Hamiltonian action and factorizable reduction}\label{double}
We will consider the action of the Manin triple $(\tmf g, \mf g, \hmf g)$
defined by a Poisson Lie group $G$ (cf. theorem \ref{app:liedouble}, also
\cite{Lu}), with dual group $\hat G$. 
In this context, the notion of Hamiltonian action on an extended K\"ahler
manifold is
\begin{defn}\label{double:biactdef}
The (infinitesimal) (left) action of $\tmf g$ on $M$ is \emph{bi-Hamiltonian} if it is
induced by a (left) $\JJC_1$-Hamiltonian action of $G$ together with a (left)
$\JJC_2$-Hamiltonian action of $\hat G$.
\end{defn}
We will use $\mu$ and $\hat \mu$ to denote the moment maps of the $G$ and $\hat
G$ actions, respectively.
Suppose that $G$ is a factorizable Poisson Lie group (definition
\ref{app:factor}). Let $S : \hat G \dto G$ be the local diffeomorphism defined
by $\tar s$ and the exponential maps at the identity elements $\hat e \in \hat
G$  and $e \in G$, then $dS(\hat e) = \tar s$.
We consider the reduction by the bi-Hamiltonian action of $\tilde G$.
\begin{assumption}\label{double:assum} We will need the following conditions:
\begin{enumerate}
\item[$(0)$] In the following, $G$ is always a factorizable Poisson Lie group.
\item The identity elements $e \in G$ and $\hat e \in \hat G$ are regular values
of $\hat \mu$ and $\mu$ respectively.
\item $\hat \mu^{-1}(e) = \mu^{-1}(\hat e)$ and is denoted $M_0$.
\item Restricted over the identity elements, $d\hat \mu = dS \circ d\mu( = \tar
s \circ d\mu)$.
\end{enumerate}
\end{assumption}
\begin{remark}\label{double:property}
\rm{ It follows from $(2)$ that $M_0$ is preserved by the $\tilde G$ action.
By $(3)$, we see that $\hat \mu^* = \mu^*\circ \tar s^* = \mu^*\circ \tar s$
when restricted to $M_0$, since $s$ is symmetric. Thus on $M_0$ we have
$$\hat \mu^* \theta_{\hat \tau} = \mu^*\circ \tar s (\theta_{\hat \tau}) = \mu^*
\hat\theta_{\tar s(\hat \tau)} \text{ for } \hat \tau \in \hmf g.$$
}
\end{remark}
\begin{lemma}\label{double:linear}
Let $\VV = V \dsum V^*$ and $(\JJ_1, \JJ_2; \GG)$ be a linear generalized
K\"ahler structure. Let $K \subset V^*$ and define
$N_j = a\circ \JJ_j(K)$ for $j = 1, 2$ where $a : \VV \to V$ is the projection.
Assume that for $j = 1, 2$
\begin{enumerate}
\item[$(1)$] $K + \JJ_j(K) \subset \Ann(K, \JJ_j(K))$,
\item[$(2)$] $\JJ_j(K) \inter V^* = \{0\}$ and
\item[$(3)$] $N_1 \inter N_2 = \{0\}$.
\end{enumerate}
We then have $\tilde \VV_K = \VV_K \dsum K$, $N_1 \dsum N_2 \subset \Ann_V(K)$
and the exact sequence
$$0 \to \frac{\Ann_{V^*}(N_1, N_2)}{K} \to \VV_K \xto{a_K} \Ann_V(K) \to 0.$$
\end{lemma}
{\it Proof: } Let $K' = \JJ_1(K) \dsum \JJ_2(K)$. It then follows from lemma
\ref{courant:missingrank}.
\qed

\begin{defn}\label{courant:extcourant}
A Courant algebroid $E$ on $M$ is an \emph{effective Courant algebroid} if it
fits in the following diagram:
\begin{equation*}
\xymatrix{
0 \ar[r] & {T^*M} \ar[r] & {E} \ar[r]^p \ar@{..>}[rd]_{a} & {E_0} \ar[r] \ar[d]&
0 \\
& & & {TM} \ar[d] \ar[rd] \\
& & & 0 & 0 \\
}
\end{equation*}
where $a$ is the anchor map and the sequences are all exact.
\end{defn}

The usual constructions of $B$-transformation for $B \in \Omega^2(M)$ and
twisting of the Courant bracket by $H \in \Omega^3_0(M)$ are valid for an
effective Courant algebroid $E$.
\begin{theorem}\label{double:courant}
Given assumption \ref{double:assum} and let $(\tilde G, G, \hat G)$ be a (local)
double Lie group whose Lie algebras form the Manin triple $(\tmf g, \mf g, \hmf
g)$, where $\tilde G$ is connected but not necessarily simply connected (compare
to theorem \ref{app:liedouble}). Suppose that the action of $\tmf g$ induces an
action of $\tilde G$, which is proper and free on $M_0$, then there is 
an effective Courant algebroid $\TTC_e \tilde Q$ 
on $\tilde Q = M_0 /\tilde G$.
\end{theorem}
{\it Proof: } Let $\Gamma(\cdot)^{\tilde G}$ denote the set of $\tilde
G$-invariant sections. By lemma \ref{poisson:equiv} we see that 
$$(\mu^*\hat\theta, \JJC_1(\mu^*\hat\theta), \JJC_2(\mu^*\hat\theta),
\GGC(\mu^*\hat\theta)) = (\mu^*\hat\theta, \JJC_1(\mu^*\hat\theta)) \dsum
\JJC_2(\mu^*\hat\theta, \JJC_1(\mu^*\hat\theta))$$
is preserved by the $G$-action. Similarly, it's also preserved by $\hat G$ and
it follows that it's preserved by the action of $\tilde G$. Analogously, the
bundles $(\mu^*\hat\theta, \JJC_1(\mu^*\hat\theta), \JJC_2(\mu^*\hat\theta))$
and $(\mu^*\hat\theta)$ are preserved by the $\tilde G$-action.
Let $\KKC = (\mu^*\hat\theta)$ and $\KKC' =
\JJC_1(\mu^*\hat\theta)\dsum\JJC_2(\mu^*\hat\theta)$, then the conditions for
lemma \ref{courant:reduce} $(2)$ are satisfied. Thus $\TTC''_\mu M_0 =
\frac{\Ann(\mu^*\hat\theta, \JJC_1(\mu^*\hat\theta),
\JJC_2(\mu^*\hat\theta))}{(\mu^*\hat \theta)}$ descends to an Courant algebroid
$\TTC_e \tilde Q$ on $Q$.

Another way to see the Courant algebroid structure is to follow theorem
\ref{kahler:reduct}. 
Using the decomposition in lemma \ref{double:linear}, we define the projection
$$\pi : \Ann(\mu^*\hat\theta, \JJC_1(\mu^*\hat\theta), \JJC_2(\mu^*\hat\theta))
\to \TTC'_\mu  M_0$$ 
and the bracket $*_\mu$:
$$\mf X *_\mu \mf Y = \pi(\mf X *_H \mf Y) \text{ for } \mf X, \mf Y \in
\Gamma(\TTC'_\mu  M_0)^{\tilde G} \subset \Gamma(\Ann(\mu^*\hat\theta,
\JJC_1(\mu^*\hat\theta), \JJC_2(\mu^*\hat\theta))^{\tilde G}.$$
By definition, $\Gamma(\TTC'_\mu  M_0)^{\tilde G}$ is closed under $*_\mu$. Then
the inclusion of $\TTC'_\mu M_0$ into $\Ann(\mu^*\hat\theta,
\JJC_1(\mu^*\hat\theta), \JJC_2(\mu^*\hat\theta))$ induces natural isomorphism
to $\TTC''_\mu  M_0$, and the brackets coincide.
\qed

\begin{corollary}\label{double:bfield}
With the same assumptions as in theorem \ref{double:courant}, let $(M, \JJC_1',
\JJC_2'; \GGC')$ be the $B_1$-transformed generalized K\"ahler structure for
$B_1 \in \Omega^2(M)^{\tilde G}$. 
Let all other choices be the same. 
Then the effective Courant algebroid $\TTC'_e \tilde Q$ induced from $(\JJC_1',
\JJC_2'; \GGC')$ is a $b$-transformation of $\TTC_e \tilde Q$, for some $b \in
\Omega^2(\tilde Q)$.
\end{corollary}
{\it Proof: }
Choose a connection form $\tilde \theta$ of the $\tilde G$-principle bundle $M_0
\to \tilde Q$ and with respect to a choice of basis of $\tmf g$ we have $\tilde
\theta_j$ and $\tilde X_j$.
Consider the form $\tilde b = \prod_j(1- \tilde\theta_j\wedge \iota_{\tilde
X_j}) B_1|_{M_0}$, where the terms in brackets are considered operators on
$\Omega^2(M_0)$. Then $\tilde b$ is horizontal with respect to $\tilde G$-action
and the transformation $e^{\tilde b}$ preserves $\Ann(\mu^*\hat\theta,
\JJC_1(\mu^*\hat\theta), \JJC_2(\mu^*\hat\theta))$. From which the result
follows.
\qed

\section{Courant and $T$-duality}\label{torus}
The \emph{Courant duality} is the following. Consider a bi-Hamiltonian action as
in definition \ref{double:biactdef}. Suppose that reduction of $G$- (resp. $\hat
G$-) action at $\hat e \in \hat G$ (resp. at $e \in G$) as given in theorem
\ref{kahler:reduct} exists and denote it $(Q, \JJC_1, \JJC_2)$ (resp. $(\hat Q,
\hat \JJC_1, \hat\JJC_2)$):
\begin{defn}\label{double:basic}
The structures $(Q, \JJC_1, \JJC_2)$ and $(\hat Q, \hat \JJC_1, \hat\JJC_2)$ are
\emph{Hamiltonian dual} to each other. When the assumptions of theorem
\ref{double:courant} holds, the structures are said to be \emph{Courant dual} to
each other.
\end{defn}
Geometrically, the Hamiltonian duality as defined above has a significant
drawback: \emph{a priori}, the level sets $M_{\hat e} = \mu^{-1}(\hat e)$ and
$M_e = \hat\mu^{-1}(e)$ might have nothing to do with each other and the
relation between the geometry and topology of the quotients $Q$ and $\hat Q$ may
not be clear. 
For Courant duality, the relation of the topology and geometry can be understood
much better.
\begin{prop}\label{torus:diagram}
Assume the conditions in theorem \ref{double:courant} and that $(\tilde G,
G, \hat G)$ is a connected double Lie group, we have the following diagram,
where the maps are principle bundles of compact Lie groups:
$$\text{
\xymatrix{
& {M_0 = Q\times_{\tilde Q} \hat Q}\ar[dl]_{\pi}^{/G} \ar[dr]^{\hat \pi}_{/\hat
G} & \\
{Q} \ar[dr]_{p}^{/\hat G} & & {\hat Q} \ar[dl]^{\hat p}_{/G}\\
& {\tilde Q} &
}
}$$
\end{prop}
{\it Proof: } 
Recall that $G \times \hat G \to \tilde G : (g, \hat g) \mapsto g \hat g^{-1}$ as well
$\hat G \times G \to \tilde G: (\hat g, g) \mapsto \hat g g^{-1}$
are diffeomorphisms for the double Lie group $(\tilde G, G, \hat G)$. The left action of $\hat G$ on $Q$ is induced from:
$$\hat g \circ g^{-1} x = \hat g g^{-1} x \text{ for } x \in M_0,$$
while the left action of $G$ on $\hat Q$ is induced from
$$g\circ \hat g^{-1} x = g\hat g^{-1} x \text{ for } x \in M_0.$$
These two actions are both free with the same quotient space $\tilde Q = M_0 / \tilde G$. 
\qed

The choice of terminology is justified by the following:
\begin{prop}\label{torus:courantiso}
With the same conditions as in proposition \ref{torus:diagram}, the Courant
algebroids on $\tilde Q$ formed by the invariant sections of $\TTC Q$ and $\TTC
\hat Q$ are isomorphic to the one defined by theorem \ref{double:courant}.
\end{prop}
{\it Proof: }
Note that the invariant sections of $\TTC Q$ lifts to $M_0$ as the $\tilde
G$-invariant sections of $\frac{\Ann(\mu^*\hat\theta, \JJC_1(\mu^*\hat
\theta))}{(\mu^*\hat\theta, \JJC_1(\mu^*\hat \theta))}$, which is isomorphic to
$\TTC''_\mu M_0$ by lemma \ref{kahler:linear}. The proposition then follows.
\qed

When the action of $G$ and $\hat G$ commute, we have:
\begin{prop}\label{torus:comute}
Let $(\tmf g, \mf g, \hmf g)$ be the Manin triple defined by a factorizable Lie
bialgebra $\mf g$. If $[\mf g, \hmf g] = 0$, then $[\tmf g, \tmf g] = 0$, i.e.
$\tmf g$ is abelian, and we write $(\tmf g, \mf g, \hmf g) = (\tmf t, \mf t,
\hmf t)$.
\end{prop}
{\it Proof: }
By \eqref{app:dsum}, we have for $\tau \in \mf g$ and $\hat \omega \in \hat g$:
$$[(\tau, \tau), (\tar{r_+}(\hat\omega), \tar{r_-}(\hat\omega))] = ([\tau,
\tar{r_+}(\hat\omega)], [\tau, \tar{r_-}(\hat\omega)]) = 0,$$
which implies that $[\tau, \tar s(\hat\omega)] = 0$. Since $\underline s$ is
invertible, we see that $\mf g$ is abelian. 
Then $(\tmf g, \mf g, \hmf g)$ form a Manin triple implies that $\tmf g$ and
$\hmf g$ are abelian as well.
\qed

Because of this, in the following we work under the assumption
\ref{double:assum} and that $\tilde G$ is abelian. The notations $\tilde T$, $T$
and $\hat T$ will mean that the respective groups are compact, i.e. torus.
\begin{lemma}\label{double:torusprop}
Both $\JJC_1$ and $\JJC_2$ are preserved by the $\tilde T$-action. For any $\tau
\in \mf t$ and $\hat \omega \in \hmf t$, we have
$d\<\JJC_1(\mu^*\hat\theta_\tau), \JJC_2(\hat\mu^*\theta_{\hat \omega})\> = 0$.
We define the pairing
\begin{equation}\label{double:geopairing}
P : \mf t \tensor \hmf t \to \RR : \tau\tensor \hat \omega \mapsto
2\<\JJC_1(\mu^*\hat\theta_\tau), \JJC_2(\hat\mu^*\theta_{\hat \omega})\>
\end{equation}
then $P$ is non-degenerate, i.e. $\tau = 0 \in \mf t \iff P(\tau, \hat \omega) =
0$ for all $\hat \omega \in \hmf t$ and vice versa for $\hat \omega$.
\end{lemma}
{\it Proof: }
By definition, the $\mf t$-action preserves $\JJC_2$ and $\hmf t$-action
preserves $\JJC_1$. Then by the proof of lemma \ref{poisson:equiv} and $d\hat
\theta = 0$ from abelian-ness, it follows that $\JJC_1$ as well as $\mu^*(\hat
\theta_\tau)$ are preserved by $\mf t$. Thus the $\tilde T$-action preserves
$\JJC_1$. Similarly, $\JJC_2$ and $\hat\mu^*(\theta_{\hat \omega})$ are
preserved by the $\tilde T$-action. 
We have:
$$\JJC_1(\mu^*\hat\theta_\tau)*_H \JJC_2(\hat\mu^*\theta_{\hat \omega}) = 0
\text{ and } \JJC_2(\hat\mu^*\theta_{\hat \omega}) *_H
\JJC_1(\mu^*\hat\theta_\tau) = 0.$$
Add the above two equations, we see that along $M_0$,
$$d\<\JJC_1(\mu^*\hat\theta_\tau), \JJC_2(\hat\mu^*\theta_{\hat \omega})\> =
0.$$
That $P$ is non-degenerate follows from non-degenerate-ness of the extended
metric $\GGC$.
\qed

As corollary of theorem \ref{kahler:reduct} and theorem \ref{double:courant}, we
note that for bi-Hamiltonian action of $(\tilde T, T, \hat T)$ where all groups
acting properly and freely, the reduction as described in theorem
\ref{double:courant} factorizes in two ways, via $Q$ or $\hat Q$, thus
\begin{defn}\label{torus:courant}
The structures $\TTC Q$ and $\TTC \hat Q$ are said to be \emph{(Courant)
$T$-dual} to each other.
\qed
\end{defn}
\begin{assumption}\label{torus:assum}
In the rest of this section, we assume that the $\tilde T$-action preserves a
splitting of $\TTC M$.
\end{assumption}
Consider the reduced structures on $Q = M_0 / T$ and $\hat Q = M_0/\hat T$. By
corollary \ref{courant:split}, 
the structures are both twisted generalized K\"ahler structures, whose twisting
form can be described with a choice of connection forms. Let $\tilde \Theta$ be
a connection form on $M_0$ as principle $\tilde T$-bundle. 
Choose basis $\{\tau_j\}$ and $\{\hat \tau_j\}$ of $\mf t$ and $\hmf t$
respectively, and denote $\theta_j, X_j + \xi_j, \Theta_j$ and $\hat \theta_j,
\hat X_j + \hat \xi_j, \hat \Theta_j$ the corresponding components. 
We define:
\begin{equation}\label{torus:btrans}
\tilde B = B + \hat B = \left(\Theta\wedge \xi -
\frac{1}{2}\sum_{j,k}\Theta_j\wedge\Theta_k \cdot \iota_{X_k}\xi_j\right) +
\left(\hat \Theta\wedge \hat \xi - \frac{1}{2}\sum_{j,k} \hat\Theta_j \wedge
\hat\Theta_k \cdot \iota_{\hat X_k}\hat \xi_j\right).
\end{equation}
Then $\tilde B$ is $\tilde T$-invariant on $M_0$. When the actions generated by
$\mf t$ and $\tmf t$ are proper, the forms $\Theta$ and $\hat\Theta$ become
connection forms on $M_0$ as respectively $T$ and $\hat T$ principle bundles.
\begin{theorem}\label{torus:duality}
Suppose assumption \ref{double:assum} holds and let $(M, \JJ_1^{\tilde B},
\JJ_2^{\tilde B}; \GG^{\tilde B})$ be the $-\tilde B$-transformed structure on
$M$, with $\tilde B$ defined above. Then the induced Courant algebroid on
$\tilde Q = M_0/\tilde T$ remains  unchanged. Let $h$ (resp. $\hat h$) be the
twisting form of the corresponding reduced structure on $Q$ (resp. $\hat Q$),
then
\begin{equation}\label{torus:dualeq}
\hat\pi^*\hat h - \pi^*h = d(\hat\Theta \wedge \Theta),
\end{equation}
where on the right hand side we use also the pairing \eqref{double:geopairing}.
\end{theorem}
{\it Proof: }
Direct computation shows that the $\tilde T$-horizontal part of $\tilde B$ is
$0$. Thus by corollary \ref{double:bfield}, the Courant algebroid structure on
$\tilde Q$ remains unchanged under the $-\tilde B$-transformation.

The $-\tilde B$-transformed structures on $M$ has twisting form $\tilde H = H +
d\tilde B$.
Let $\JJ_1^{\tilde B} = e^{-\tilde B} \JJ_1 e^{\tilde B}$ and so on. We compute
\begin{equation*}
\iota_{X_l}\tilde B  = \xi_l - \hat \Theta \cdot \iota_{X_l}\hat \xi = \xi_l -
\sum_j \hat \Theta_j \cdot \iota_{X_l} \hat \xi_j
\end{equation*}
and it follows that $\JJ_1^{\tilde B}(\mu^*\hat \theta_l) = X_l + \xi_l'$ where
$\xi_l' = \sum_j \hat \Theta_j \cdot \iota_{X_l} \hat \xi_j$. We note that
$\iota_{X_j}\xi_l' = 0$, and the twisting form $h$ satisfies
$$\pi^*h = \tilde H + d(\Theta \wedge\xi').$$
Similarly, we have $\hat \pi^*\hat h = \tilde H + d(\hat \Theta \wedge \hat
\xi')$ where $\hat \xi_l' = \sum_j \Theta_j \cdot \iota_{\hat X_l} \xi_j$.
More explicitly, we compute
\begin{equation*}
\begin{split}
\hat\pi^*\hat h - \pi^*h = & d\left(\sum_{j,k}\hat\Theta_j \wedge \Theta_k \cdot
\iota_{\hat X_j}\xi_k - \sum_{j,k} \Theta_k \wedge \hat\Theta_j \cdot
\iota_{X_k} \hat\xi_j\right) \\
= & d\sum_{j,k}\hat \Theta_j \wedge \Theta_k \cdot (\iota_{\hat X_j} \xi_k +
\iota_{X_k} \hat \xi_j) \\
= & d\sum_{j,k}\hat \Theta_j \wedge \Theta_k \cdot 2\<\JJ_2(\hat \mu^*\theta_j),
\JJ_1(\mu^*\hat \theta_k)\> \\
= & d(\hat\Theta \wedge \Theta).
\end{split}
\end{equation*}
where the last step we use the pairing $P$ as given in
\eqref{double:geopairing}.
\qed
\begin{remark}\label{torus:dualtrans}
\rm{ We note that the $T$- or $\hat T$-horizontal part of $\tilde B$ in general
do not vanish. Thus the structures on $Q$ and $\hat Q$ are $B$-transformed from
their respective original structures. 
With proposition \ref{torus:courantiso} the theorem above states that the
Courant algebroid on $\tilde Q$ formed by the set of invariant sections of $\TTC
Q$ or $\TTC \hat Q$ are still isomorphic to the original one. We note also that
the equation \eqref{torus:dualeq} coincides with the equation in the physics
literature, where $M_0$ is to be the correspondence space of the $T$-dual
bundles $Q$ and $\hat Q$. 
It's shown (e.g. \cite{Bouwknegt}) that the twisted cohomology of $T$-dual
principle bundles are isomorphic. Since the twisted cohomology only depends on
the cohomology class of the twisting, the same is true for the structures on $Q$
and $\hat Q$ before applying $B$-transformation. 
In \cite{Cavalcanti}, proposition \ref{torus:courantiso} is shown when $Q$
and $\hat Q$ are $T$-dual $S^1$-principle bundles, with twisted generalized
complex structures, by directly defining the isomorphism.
}
\end{remark}
\subsection{Example}\label{torus:example}
The following example is considered in \cite{Hu} and we recall the setup and
point out its relevance to the current discussion. Let $M = \CC^2
\setminus\{(0,0)\}$ and consider the coordinates $z = (z_1, z_2) = (x_1 + i y_1,
x_2 + iy_2) = (x_1, y_1, x_2, y_2)$. Let $r^2 = |z_1|^2 + |z_2|^2$ and  $J =
\left(\begin{smallmatrix}0 & -1 \\ 1 & 0 \end{smallmatrix}\right)$. Consider the
following structures:
\begin{equation*}
\JJ_1 =  \(\begin{matrix}
	0 & 0 & r^2 J & 0 \\
	0 & -J & 0 & 0 \\
	r^{-2} J & 0 & 0 & 0 \\
	0 & 0 & 0 & -J
	\end{matrix}\),
\JJ_2 = \(\begin{matrix}
	J & 0 & 0 & 0 \\
	0 & 0 & 0 & -r^2 J \\
	0 & 0 & J & 0 \\
	0 & -r^{-2} J & 0 & 0
	\end{matrix}\),
\end{equation*}
where the labelling on rows are $(T z_1, T z_2, T^*z_1, T^*z_2)^T$. Then $(M,
\JJ_1, \JJ_2)$ is an $H$-twisted generalized K\"ahler structure where
$$H = -\sin (2\lambda) d\lambda \wedge d\phi_1 \wedge d\phi_2$$
in the polar coordinates $(z_1, z_2) = r(e^{i\phi_1}\sin \lambda, e^{i\phi_2}
\cos \lambda)$. In particular, $[H] \neq 0 \in H^3(M)$ (cf. \cite{Gualtieri}).

Let $\tilde T = S^1 \times S^1$ and $(e^{i\theta_1}, e^{i\theta_2})$ be the
coordinates. It acts on $M$ via 
$$(e^{i\theta_1}, e^{i\theta_2})\circ (z_1, z_2) = (e^{i\theta_1} z_1,
e^{-i\theta_2} z_2).$$
Let $T$ and $\hat T$ be the first and second $S^1$ respectively, then $(\tilde
T, T, \hat T)$ is a double Lie group and the action of $\tilde T$ is
bi-Hamiltonian and satisfies assumption \ref{double:assum} with the common
moment map $f = \ln r$. The $T$ and $\hat T$ actions are generated respectively
by 
$$\JJ_1(df) = \frac{\partial}{\partial \phi_1} - \cos^2 \lambda d\phi_2 \text{
and } \JJ_2(df) = -\frac{\partial}{\partial \phi_2} + \sin^2 \lambda d\phi_1.$$
It follows that $2\<\JJ_1(df), \JJ_2(df)\> = 1$, i.e. $T$ and $\hat T$ are dual
tori in the standard sense. We note that the actions are not free. Consider the
submanifold $M' = M\setminus \left(\{z_1 = 0\} \union \{z_2 = 0\}\right)$, on
which the actions are free. The reduced structures of the $T$ and $\hat T$
action on $M'$ are respectively the opposite and standard K\"ahler structures on
$D^2 \setminus \{0\}$. 
Our results then state that they are $T$-dual to each other.

\section{$T$-duality group}\label{dual}
We first consider the linear case and use the notations and assumptions of lemma
\ref{double:linear}. We note that the natural pairing on $\VV$ induces a pairing
$P_K$ on $\JJ_1(K)\dsum \JJ_2(K)$, which can also be seen as induced from the
pairing $\<\cdot, \GG(\cdot)\>$ defined on $K$ by $\GG$. By the positive
definiteness of $\GG$ we see that $P_K$ has signature $(m,m)$ where $m = \dim
K$. Completely parallel to lemma \ref{kahler:sequence}, we have
\begin{lemma}\label{dual:linear}
We use the notations and assumptions of lemma \ref{double:linear}. Let $K'
\subset \JJ_1(K)\dsum \JJ_2(K)$ be a maximal isotropic subspace with respect to
$P_K$ and $N' = a(K')$ where $a: \WW \to V$ is the projection, then there is a
self-dual exact sequence:
$$0 \to W^*_{K'} \xto{a_{K'}^*} \WW_K \xto{a_{K'}} W_{K'} \to 0, \text{ with }
W_{K'} =\frac{\Ann_V(K)}{N'} \text{ and } W^*_{K'} = \frac{\Ann_{V^*}(N')}{K},$$
where $\WW_K = \Ann(K, K', \GG(K'), \GG(K)) = \Ann(K, \JJ_1(K), \JJ_2(K),
\GG(K))$. 
\end{lemma}
{\it Proof: }
By the conditions in lemma \ref{double:linear}, we see that $K'\inter V^* =
\{0\}$, $K\dsum K' \subset \Ann(K, K')$. Let $a_{K'}$ be the map induced from
$a$. Let $\UU_{K'} = \Ann(K, K')$, then the kernel of the induced map $\UU_{K'}
\to W_{K'}$ is $\Ann_{V^*}(N') \dsum K'$ and thus the kernel of $a_{K'}$ is
$$\ker a_{K'} = (\Ann_{V^*}(N') \dsum K') \inter \VV_K.$$
By the decomposition $\UU_{K'} = \WW_K \dsum (K\dsum K')$ and inclusion $K \dsum
K' \subset \Ann_{V^*}(N') \dsum K'$, we see that $\ker a_{K'} \simeq
\frac{\Ann_{V^*}(N')}{K}$. Since $\ker a_{K'}$ is maximally isotropic with
respect to the induced pairing $\<,\>_K$ on $\WW_K$, we see that the exact
sequence is self-dual.
\qed

Using the notations in (the proof of) theorem \ref{kahler:reduct}, we have
\begin{prop}\label{dual:kahler}
Under the condition of theorem \ref{double:courant} and let $T' \subset \tilde
T$ be a maximally isotropic subtorus of $\tilde T$ with respect to the pairing
$P$ as in \eqref{double:geopairing}, i.e., the Lie algebra $\mf t'$ is a
Lagrangian subspace of $\tmf t$. Then the reduced space $Q' = M_0 /T'$ has a
natural extended K\"ahler structure. 
\end{prop}
{\it Proof: }
By the proof of theorem \ref{double:courant}, the bundles $(\mu^*\hat \theta)$,
$(\mu^*\hat \theta, \JJC_1(\mu^*\hat \theta), \JJC_2(\mu^*\hat\theta))$ and
$\TTC_\mu M_0$ are all preserved by the $\tilde T$-action, and thus are
preserved by the $T'$-action. Let $\KKC = (\mu^*\hat\theta)$ and $\KKC'$ be the
subbundle generated by the infinitesimal fields $\{X_{\tau'} + \xi_{\tau'} |
\tau' \in K'\}$, then it follows from the proofs of lemmata \ref{poisson:moment}
and \ref{double:torusprop} that $(\mu^*\hat \theta, \KKC')$ is preserved by the
$T'$-action. Since $T'$ is isotropic, we have $\KKC \dsum \KKC' \subset
\Ann(\KKC, \KKC')$ and lemma \ref{courant:reduce} $(1)$ gives an extended
tangent bundle $\TTC Q'$ on $Q'$.
It follows from proposition \ref{double:torusprop} that $\JJC_1$ and $\JJC_2$
are both invariant with respect to the $T'$-action and thus descend to $\JJC'_1$
and $\JJC'_2$ on $\TTC Q'$, which define an extended K\"ahler structure.
\qed

The group $O(m,m;\ZZ)$ is called the \emph{$T$-duality group} in the physics
literature (\cite{Mathai1}, \cite{Hull} and the references therein). In the
physics literature, for each element of $O(m,m;\ZZ)$ it is associated a pair of
related $T$-dual principle bundles with $H$-fluxes. The physical theory on such
related structures are expected to be the same. In our construction, the
following holds.
\begin{corollary}\label{dual:tgroup}
Suppose that the action of $\tilde T$ preserves a splitting of $\TTC M$. Let $g
\in O(m,m;\ZZ)$ and consider the pair of Lagrangian subgroups $T_g$ and $\hat
T_g$ with Lie algebra $g(\mf t)$ and $g(\hmf t)$. Let $Q_g$ and $\hat Q_g$ be
the reduction of $M_0$ by the groups $T_g$ and $\hat T_g$ respectively. Then the
induced structures on $Q_g$ and $\hat Q_g$ are twisted generalized K\"ahler
structures and the equation \eqref{torus:dualeq} holds for this pair after
applying certain $B$-transformation. The Courant algebroid on $\tilde Q$ defined
by the $\hat T_g$-invariant sections of $\TT Q_g$ is isomorphic to the one given
by theorem \ref{double:courant}.
\end{corollary}
{\it Proof: } Similar to \eqref{torus:btrans}, we choose basis $\{\tau^g_j\}$
and $\{\hat\tau^g_j\}$ of $\mf t^g$ and $\hmf t^g$ respectively and let
$X^g_j+\xi^g_j$, $\Theta^g_j$ and $\hat X^g_j+\hat\xi^g_j$, $\hat\Theta^g_j$ be
the corresponding components, and define
\begin{equation*}
\tilde B^g = B^g + \hat B^g = \left(\Theta^g\wedge \xi^g -
\frac{1}{2}\sum_{j,k}\Theta^g_j\wedge\Theta^g_k \cdot
\iota_{X^g_k}\xi^g_j\right) + \left(\hat \Theta^g\wedge \hat \xi^g -
\frac{1}{2}\sum_{j,k} \hat\Theta^g_j \wedge \hat\Theta^g_k \cdot \iota_{\hat
X^g_k}\hat \xi^g_j\right).
\end{equation*}
In particular, the basis $\{\tau^g_j\}$ and $\{\hat\tau^g_j\}$ can be taken as
the transformation of $\{\tau_j\}$ and $\{\hat\tau_j\}$ by $g$. The proof of
\eqref{torus:dualeq} is then completely parallel to that of theorem
\ref{torus:duality}. The isomorphism of courant algeboids is straight forward.
\qed

\subsection{Example} \label{dual:bfield}
We consider in detail the special case when $g = e^{b}$ where $b : \mf t \to
\hmf t$ is skew-symmetric with respect to the pairing $P$. Then $g(\mf t) =
graph(b)$ and $g(\hmf t) = \hmf t$. Let $\{\tau_j\}$ and $\{\hat \tau_j\}$ be
basis of $\mf t$ and $\hmf t$ respectively and $(b_{ij})$ the matrix of $b$ with
respect to these basis. Then $\{\tau^b_j = \tau_j + \sum_{k} b_{jk} \hat
\tau_k\}$ is a basis of $g(\mf t)$, where $b(\tau_j) = \sum_{k} b_{kj} \hat
\tau_k \in \hmf t$. Let $^b$ to denote the objects for the transformed
structures, then
\begin{equation*}
\left\{\begin{split}
\Theta^b_j & = \Theta_j \\
\hat\Theta^b_j & = \hat\Theta_j - \sum_k b_{jk} \Theta_k
\end{split}
\right.
\text{ and }
\left\{
\begin{split}
\hat X^b_j + \hat\xi^b_j & = \hat X_j + \hat\xi_j\\
X^b_j + \xi^b_j & = X_j + \xi_j + \sum_k b_{kj}(\hat X_k + \hat \xi_k).
\end{split}
\right.
\end{equation*}
Direct computation gives
\begin{equation*}
\left. \begin{split}
(\pi^b)^*h^b = & \pi^*h - d\sum_{j,k,l}b_{lj}\Theta_k \wedge \Theta_j
\<X_k+\xi_k, \hat X_l + \hat\xi_l\> \\
(\hat \pi^b)^*\hat h^b = & \hat\pi^*\hat h - d\sum_{j,k,l}
b_{lj}\Theta_j\wedge\Theta_k \<X_k+\xi_k, \hat X_l + \hat\xi_l\>\\
\end{split} \,\,\,\, \right\} \text{ which implies }
\end{equation*}
$$(\hat \pi^b)^*\hat h^b - (\pi^b)^*h^b = d\sum_{j,k}\hat\Theta^b_j \wedge
\Theta^b_k \cdot (\iota_{\hat X^b_j} \xi^b_k + \iota_{X^b_k} \hat\xi^b_j),$$
i.e., the equation \eqref{torus:dualeq} holds for the pair of reduced structures
$Q^b$ and $\hat Q^b$. Since $\hmf t^b = \hmf t$, we have $\hat Q^b = \hat Q$,
while the twisting form $\hat h$ is changed by an exact term. As the situation
for $\mf t$ and $\hmf t$ is symmetric, we may consider $e^\beta$ for
skew-symmetric $\beta: \hmf t \to \mf t$ and obtain similar result.

\subsection{Example} \label{dual:example}
The example discussed in \S\ref{torus:example} does not admit interesting
$T$-duality group action, because $O(1, 1; \ZZ) = \{\pm 1,
\pm\left(\begin{smallmatrix}0 & 1 \\ 1 & 0 \end{smallmatrix}\right)\}$. This can
be compensated by considering a product of these, e.g. twisted structures on
$M^2$, for example, and apply \S\ref{dual:bfield}. Instead, here we consider
another situation which is not quite covered by $T$-duality group.
For the example in \S\ref{torus:example} we consider the anti-diagonal action
generated by 
$$X_d + \xi_d = \left(\frac{\partial}{\partial \phi_1}+ \frac{\partial}{\partial
\phi_2}\right) - (\cos^2\lambda d\phi_2 + \sin^2 \lambda d\phi_1),$$
then $\iota_{X_d}\xi_d = -1$. Let $\KKC = (d\mu)$ and $\KKC' = (X_d +\xi_d)$,
then the condition for lemma \ref{courant:reduce} $(1)$ does not hold. On the
other hand, 
the condition for lemma \ref{courant:reduce} $(2)$ holds and there is an induced
effective Courant algebroid 
on the corresponding reduced space, i.e. $S^2$. A more general result holds:
\begin{prop}\label{dual:nondegen}
Let $T^+ \subset \tilde T$ be a \emph{non-degenerate} subtorus of $\tilde T$
with respect to $P$, i.e. the restriction of $P$ on its Lie algebra $\mf t^+$ is
non-degenerate, then there is a natural effective Courant algebroid on the
reduced space $Q^+ = M_0 / T^+$.
\qed
\end{prop}

\section{Appendix A: Reduction of extended tangent bundle}\label{courant}
Special case of the reduction of Courant algebroid has been discussed implicitly
in our paper \cite{Hu} in showing that extended complex structure exists as the
result of reduction of generalized complex manifold and in general it has been
discussed explicitly in the works \cite{Bursztyn}, \cite{Stienon} and \cite{Vaisman}. For the sake
of completeness, we prove the reduction of Courant algebroid in the relevant
context of our construction in this article, i.e. for extended tangent bundles. 
We will use the notations in \S\ref{recall}.
\begin{lemma}\label{courant:extend}
Let $\VV = V \dsum V^*$ with the natural pairing $\<,\>$, $K \subset V^*$ and
$K' \subset \VV$ so that $K' \inter V^* = \{0\}$. Define
$N' = a\circ K'$ where $a : \VV \to V$ is the projection. Assume that
\begin{enumerate}
\item[$(1)$] $K + K' \subset \Ann(K, K')$ and
\item[$(2)$] $K' \inter V^* = \{0\}$.
\end{enumerate}
Let $\VV_K = \frac{\Ann(K, K')}{(K, K')}$, then $\<,\>$ descends to
non-degenerate pairing $\<,\>_K$ on $\VV_K$ and we have the self-dual exact
sequence
$$0 \to W_K^* \to \VV_K \xto{a_K} W_K \to 0 \text{ where } W_K =
\frac{\Ann_V(K)}{N'}.$$
\end{lemma}
{\it Proof: } See \cite{Hu}, lemma $4.3$.
\qed
\begin{lemma}\label{courant:missingrank}
Using the same notations as in lemma \ref{courant:extend} and replacing
assumption $(1)$ by one of the following statements that are equivalent to each
other: 
\begin{enumerate}
\item[$(1')$] $K \subset \Ann(K, K')$ and $\<,\>$ induces a non-degenerate
pairing on $\VV_K' = \frac{\Ann(K,K')}{K}$,
\item[$(1'')$] $K \subset \Ann(K, K')$ and $\<,\>$ restricts to a non-degenerate
pairing on $K'$.
\end{enumerate}
then we have the exact sequence:
$$0 \to W^*_K \to \VV_K' \xto{a_K'} \Ann_V(K) \to 0.$$
\end{lemma}
{\it Proof: } The surjectivity of $\Ann(K, K') \to \Ann_V(K)$ is easy and
everything then follows.
\qed

\begin{defn}\label{courant:idea}
Let $S$ be a subspace of sections in $\TTC M$ which is closed with respect to
$*$. A closed subspace $S' \subset S$ is a \emph{two-sided null ideal} if $S *
S' \subset S'$, $S' * S \subset S'$ and $\<S',S\> = 0$.
\end{defn}
It follows that when $S'$ is a two-sided null idea of $S$, the structure $(S, *,
\<,\>)$ induces one such structure on the quotient space $S/S'$, which also
satisfies \eqref{app:courantdefn2} and \eqref{app:courantdefn3}.
\begin{lemma}\label{courant:reduce}
Let $(M, \TTC M)$ be a manifold with an extended tangent bundle $\TTC M$ and
$M_0 \subset M$ a submanifold.
Let $\KKC \subset T^*M|_{M_0}$ and $\KKC' \subset \TTC M|_{M_0}$ be two
subbundles of rank $m$ and $m'$ respectively so that $T M_0 = \Ann_{T M}(\KKC)$
and $\KKC' \inter T^*M = \{0\}$. Suppose that $\KKC$ is generated by sections
$\{\theta_j\}_{j = 1}^m$ so that $d\theta_j \in \Gamma(\wedge^2 \KKC)$ and
$\KKC'$ is generated by sections $\{\mf X_j\}_{j = 1}^{m'}$. Let $\tilde\sigma$
be the infinitesimal action generated by $\{\mf X_j\}_{j = 1}^{m'}$ via the
Loday bracket $*$. Suppose that $\Ann(\KKC, \KKC')$ is preserved by
$\tilde\sigma$.

If furthermore, we suppose that the action $\tilde\sigma$ on $M_0$ is induced by
a morphism $G \to \tilde\GGS_H$ where $G$ is compact of dimension $m'$ and the
geometrical action $\sigma$ is free. Let $Q = M_0 / G$ then
\begin{enumerate}
\item If $(\KKC, \KKC')\subset \Ann(\KKC, \KKC')$, then $\frac{\Ann(\KKC,
\KKC')}{(\KKC, \KKC')}$ descends to an extended tangent bundle on $Q$.
\item If $\KKC \subset \Ann(\KKC, \KKC')$ is preserved by $\sigma$ and $\<,\>$
induces a non-degenerate pairing on $\KKC'$, then $\frac{\Ann(\KKC,
\KKC')}{\KKC}$ descends to an effective Courant algebroid on $Q$.
\end{enumerate}
\end{lemma}
{\it Proof: } Let $a : \TTC M \to TM$ be the projection. Let $\mf X, \mf X'
\in \Gamma(\Ann(\KKC, \KKC'))$, then
$$ \<\mf X_j, \mf X\> = 0 \text{ and } \iota_X \theta_j = 0 .$$ 
From the assumptions, $\Ann_{TM}(\KKC)$ is an integrable distribution and $M_0$
is a leaf of this distribution. 
We also have $\<\mf X *_H \mf X', \theta_j\> = \iota_{[X, X']} \theta_j = 0$ by
assumption on $\theta_j$. It then follows that
$$\<\mf X *_H \theta_j, \mf X'\> = X'\<\mf X, \theta_j\> - \<\mf X*_H \mf X',
\theta_j\> = 0,$$
i.e. $\mf X*_H \theta_j \in \Gamma(\KKC, \KKC')$. Similarly we have $ \theta_j
*_H \mf X \in \Gamma(\KKC, \KKC')$.
Let $\NNC' = a(\KKC')$, then $\NNC' \subset \Ann_{TM}(\KKC)$. We see that the
(geometrical) action of $G$ preserves $M_0$ and the quotient $Q$ is
well-defined.

\vspace{0.1in}
\noindent
{\bf Case $\mathbf{(1)}$.} Let 
$$S_1 = \{\mf X \in \Gamma(\Ann(\KKC, \KKC')) | \mf Y *_H \mf X \in \Gamma(\KKC,
\KKC') \text{ for all } \mf Y \in \Gamma(\KKC')\}.$$ 
For any $\mf X, \mf X' \in S_1$, $\mf Y \in \Gamma(\KKC')$ and $\mf Z \in
\Gamma(\Ann(\KKC, \KKC'))$
we compute 
\begin{equation*}\begin{split}
& \<\mf X *_H \mf X', \mf Y\> = a(\mf X) \<\mf X', \mf Y\> - \<\mf X', \mf X *_H
\mf Y\> = \<\mf X', \mf Y *_H \mf X\> - a(\mf X')\<\mf Y, \mf X'\> = 0,\\
& \<\mf X *_H \mf Y, \mf Z\> = a(\mf Z) \<\mf X, \mf Y\> - \<\mf Y *_H \mf X,
\mf Z\> = 0. \\
\end{split}
\end{equation*}
Thus we have $\mf X *_H \mf X' \in \Gamma(\Ann(\KKC, \KKC'))$ as well as $\mf X
*_H \mf W$ and $\mf W *_H \mf X \in \Gamma(\KKC, \KKC')$ for all $\mf W \in
\Gamma(\KKC, \KKC')$. It follows that
$$\mf Y *_H (\mf X *_H \mf X') = (\mf Y *_H \mf X)*_H \mf X' + \mf X *_H (\mf Y
*_H \mf X') \Rightarrow \mf X *_H \mf X' \in S_1,$$
i.e. $S_1$ is closed under $*_H$ and $\Gamma(\KKC, \KKC')$ is a two-sided null
ideal in $S_1$.
Thus the structures $(*_H, \<,\>)$ descend to $(*^G, \<,\>^G)$ on
$\frac{S_1}{\Gamma(\KKC, \KKC')} = \Gamma(\frac{\Ann(\KKC, \KKC')}{(\KKC,
\KKC')})^G$. Lemma \ref{courant:extend} implies that $\frac{\Ann(\KKC,
\KKC')}{(\KKC, \KKC')}$ descends to an extension $E_Q$ of $TQ$ by $T^*Q$. The
equations \eqref{app:courantdefn1}, \eqref{app:courantdefn2} and
\eqref{app:courantdefn3} holds by the comment above.

\vspace{0.1in}
\noindent
{\bf Case $\mathbf{(2)}$.} Let
$$S_2 = \{\mf X \in \Gamma(\Ann(\KKC, \KKC')) | \mf X_j *_H \mf X \in
\Gamma(\KKC)\}.$$
Completely parallel to case $(1)$ above, we see that $S_2$ is closed under
$*_H$. By definition $\Gamma(\KKC)$ is a two-sided null ideal in $S_2$ with
respect to $(*_H, \<,\>)$. 
By lemma \ref{courant:missingrank}, $T^*Q \subset \ker a'_Q$. 
\qed

Let $\Theta$ be a connection form on $\pi : M_0 \to Q$, and $\Theta_j$
the component dual to $X_j$.
\begin{corollary}\label{courant:split}
In the above lemma, suppose that the action of $G$ on $\TTC M$ preserves a
splitting into $\TT M$ with $H$-twisted structures. Let $\mf X_j = X_j + \xi_j$
under the splitting. Then $\TTC Q$ splits into $\TT Q$ with $h$-twisted
structures, where $\pi^* h = H + dB$ with 
$$B = \Theta\wedge \xi_\mu - \frac{1}{2}\sum_{j,k}\Theta_j\wedge \Theta_k \cdot
\iota_{X_k}\xi_j.$$
\end{corollary}
{\it Proof:} The action preserving the splitting implies that $d\xi_j =
\iota_{X_j} H$ and $\LLC_{X_j} H = 0$. It follows that $\LLC_{X_\tau} B = 0$ for
all $\tau \in \mf g$. Let $B' = \prod_j(1 - \Theta_j \wedge \iota_{X_j}) B$ be
the horizontal part of $B$, where $\Theta_j \wedge \iota_{X_j}$ is interpreted
as an operator on $\Omega^2(M)$. Direct computation gives $\iota_{X_\tau} B =
\xi_\tau$ and $B' = 0$. Apply $B$-transformation (or choose a different
splitting), we have $\Ann(\KKC, \KKC') \mapsto \Ann(\KKC, \{X_j\})$ and it
defines a splitting of $\TTC Q$. Under the $B$-transformation, the twisting form
becomes $H' = H + dB$ and we compute 
$$\iota_{X_j}(H + dB)  = \iota_{X_j}H + \LLC_{X_j}B - d\iota_{X_j}B =
\iota_{X_j} H - d\xi_j = 0.$$
It follows that there is $h \in \Omega^3_0(Q)$, so that $\pi^* h = H + dB$,
which gives the twisting form of the induced splitting of $\TTC Q$.
\qed

\begin{remark}\label{poisson:gencplxrmk}
\rm{ We note that $\LLC_{X_j} H = 0$ and $d\xi_j = \iota_{X_j} H$ for all $j$
implies that $d_G(H + \sum_{j}\xi_j u_j) = 0$, where $d_G$ is the equivariant
differential in the equivariant Cartan complex. Then $h$ in the above gives an
explicit description of the image of $[H + \sum_{j}\xi_j u_j]$ under the
isomorphism $H_G(M_0) \xto{\simeq} H(Q)$. From here, it again follows that $[h]$
is independent of the choice of the connection form.
}
\end{remark}

\section{Appendix B : Poisson Lie group and actions} \label{app:plie}
The material in this subsection is taken from \cite{Lu, Lu1} and \cite{Chari} 
(the first three chapters). More details can be found there as well as the 
references therein. We follow the convention \ref{app:conv}.
\begin{defn}\label{app:pliedef}
A Lie group $G$ is called a \emph{Poisson Lie group} if it is also a Poisson
manifold such that the multiplication map $m : G\times G \to G$ is a Poisson
map, where $G\times G$ is equipped with the product Poisson structure.
\end{defn}
Let $\pi$ be a multiplicative Poisson tensor on $G$, then $\pi|_e = 0$ where $e
\in G$ is the identity, and the linearization of $\pi$ at $e$ defines on $\hmf g
= \mf g^*$ a structure of Lie algebra $[,]^\roof$. From \cite{Weinstein},
\begin{theorem}\label{app:invariantalg}
The right (left) invariant $1$-forms on a Poisson Lie group $(G, \pi)$ form a
Lie subalgebra of $\Omega^1(G)$ with respect to the bracket
\begin{equation}\label{app:liebracform}
[\theta^\cdot_{\hat \tau}, \theta^\cdot_{\hat \omega}]^* =
-d\pi(\theta^\cdot_{\hat \tau}, \theta^\cdot_{\hat \omega}) +
\LLC_{\pi(\theta^\cdot_{\hat\tau})} \theta^\cdot_{\hat \omega} -
\LLC_{\pi(\theta^\cdot_{\hat \omega})} \theta^\cdot_{\hat \tau} =
\LLC_{\pi(\theta^\cdot_{\hat \tau})}\theta^\cdot_{\hat \omega} -
\iota_{\pi(\theta^\cdot_{\hat \omega})} d\theta^\cdot_{\hat \tau} \text{ for }
\hat\tau, \hat\omega \in \hmf g.
\end{equation}
The corresponding Lie algebra structure on $\hmf g$ coincides with the one given
by linearizing $\pi$ at the identity $e \in G$. In particular,
$\theta^\cdot_{[\hat\tau, \hat\omega]^\roof} = [\theta^\cdot_{\hat\tau},
\theta^\cdot_{\hat\omega}]^*$ for $\hat\tau, \hat\omega \in \hmf g$ and $\cdot =
l$ or $r$.
\qed
\end{theorem}
The tangent Lie algebra $\mf g$ of a Poisson Lie group $G$ is an example of Lie
bialgebra, as defined below.
\begin{defn}\label{app:liebidef}
A \emph{Lie bialgebra} is a vector space $\mf g$ with a Lie algebra structure
and a Lie coalgebra structure, which are compatible in the following sense: the
cocommutator mapping $\delta: \mf g \to \mf g \tensor \mf g$ must be a
$1$-cocycle ($\mf g$ acts on $\mf g\tensor \mf g$ by means of the adjoint
representation). A triple of Lie algebras $(\mf p, \mf p_1, \mf p_2)$ is called
\emph{Manin triple} if $\mf p$ has a nondegenerate invariant pairing $\<,\>$ and
isotropic Lie subalgebras $\mf p_1$ and $\mf p_2$ such that as vector space $\mf
p = \mf p_1 \dsum \mf p_2$.
\end{defn}
The cocommutator $\delta$ induces a Lie bracket on the dual $\hmf g$ of $\mf g$
and $(\tmf g = \mf g \dsum \hmf g, \mf g, \hmf g)$ with the natural pairing
between $\mf g$ and $\hmf g$ form a Manin triple. Conversely, the $\mf p_i$ are
dual to each other via the nondegenerate pairing $\<,\>$.
\begin{defn}\label{app:dualgp}
Let $\hat G$ be a Lie group with Lie algebra $(\hmf g, [,]^\roof)$ with a
Poisson Lie structure $\hat \pi$ so that the linearization of $\hat \pi$ at the
$\hat e \in \hat G$ gives $(\mf g, [,])$, then $(\hat G, \hat\pi)$ is \emph{a
dual} Poisson Lie group of $G$. When $\hat G$ is simply connected, the structure
$\hat \pi$ always exists and $\hat G$ is called \emph{the dual} group.
\end{defn}
\begin{defn}\label{app:doublegp}
Three Lie groups $(\tilde G; G_+, G_-)$ form a \emph{double Lie group} if
$G_{\pm}$ are both closed Lie subgroups of $G$ such that the map $G_+\times G_-
\to \tilde G : (g_+, g_-) \mapsto g_+g_-$ is a diffeomorphism. They form a
\emph{local double Lie group} if there exist Lie subgroups $G'_\pm$ of $\tilde
G$ such that $G'_i$ is locally isomorphic to $G_i$ for $i = +, -$ and the map
$G'_+\times G'_- \to \tilde G : (g'_+, g'_-) \mapsto g'_+g'_-$ is a local
diffeomorphism near identities.
\end{defn}
\begin{theorem}\label{app:liedouble}
Let $G$ be a Poisson Lie group with dual group $\hat G$, then $\mf g$ is
naturally a Lie bialgebra. Let $\tilde G$ be the connected and simply connected
Lie group with Lie algebra $\tmf g = \mf g \dsum \hmf g$ as given above, then
$(\tilde G, G, \hat G)$ form a local double Lie group.
\qed
\end{theorem}
The local double Lie group $(\tilde G, G, \hat G)$ in the theorem will be called
\emph{the local double group} of $G$. In general, if the Lie algebras of a
(local) double Lie group $(\tilde G, G, \hat G)$ coincide with the Manin triple
defined by the Lie bialgebra $\mf g$, then we say that $(\tilde G, G, \hat G)$
is \emph{a (local) double group} of $G$.

\subsection{} Let $r = \sum_i a_i \tensor b_i \in \mf g \tensor \mf g$, then it
defines a cocommutator $\delta$ via
\begin{equation}\label{app:cocommu}
\delta : \mf g \to \mf g \tensor \mf g : X \mapsto \ad_X r,
\end{equation}
which is a $1$-cocycle because it is in fact a $1$-coboundary. We write $r = s +
a$ where $s$ (resp. $a$) is the symmetric (resp. antisymmetric) part of $r$,
then 
$\delta$ as given in \eqref{app:cocommu} defines a Lie bialgebra iff
\begin{enumerate}
\item 
$s$ is $\ad$-invariant and
\item $\ldb r, r \rdb$ is $\ad$-invariant, where 
\begin{equation}\label{app:yangbaxter}
\ldb r, r \rdb  = \sum_{i,j} \left([a_i, a_j] \tensor b_i \tensor b_j + a_i
\tensor [b_i, a_j] \tensor b_j +  a_i \tensor a_j \tensor [b_i, b_j]\right).
\end{equation}
\end{enumerate}
\begin{defn}\label{app:factor} The Lie bialgebra defined by $r \in \mf g\tensor
\mf g$ as above is called a \emph{coboundary Lie bialgebra}. It is
\emph{factorizable} if $\ldb r, r \rdb = 0$ and $s$ is invertible. In this case,
$r$ is also called a \emph{factorizable $r$-matrix}.
A (local) double Lie group $(\tilde G, G, \hat G)$ is called \emph{factorizable}
if the corresponding Lie bialgebra is factorizable. In this case, we will also
call $G$ a \emph{factorizable} Poisson Lie group.
\end{defn}

For an element $r \in \mf g \tensor \mf g$, let $\tar r : \hmf g \to \mf g$ be
the map defined by $\tar r(\tau^*)(\omega^*) = (\tau^*\tensor\omega^*)(r)$.
Suppose that $r$ is factorizable and 
let $(\tmf g, \mf g, \hmf g)$ be the associated Manin triple,
then $\tmf g \simeq \mf g\dsum \mf g$ as Lie algebra.
The isomorphism is given by $\mf g \into \mf g \dsum \mf g : \tau \mapsto (\tau,
\tau)$ and
\begin{equation}\label{app:dsum}
\hmf g \into \mf g \dsum \mf g : \hat \omega \mapsto (\tar {r_+}(\hat \omega),
\tar {r_-}(\hat\omega)), \text{ with } r_\pm = a \pm s.
\end{equation}

\subsection{} It's a general fact for Poisson manifolds that
$\pi([\theta^\cdot_{\hat \tau}, \theta^\cdot_{\hat \omega}]^*) =
[\pi(\theta^\cdot_{\hat \tau}), \pi(\theta^\cdot_{\hat \omega})]$, in the
convention \ref{app:conv}. It follows that the map 
\begin{equation}\label{app:liehomo}
\rho^\cdot : \hmf g \to \Gamma(TG) : \hat\tau \mapsto \XXC^\cdot_{\hat\tau} =
\pi(\theta^\cdot_{\hat \tau})
\end{equation}
is a Lie algebra homomorphism, where $\cdot$ stands for $l$ or $r$.
\begin{defn}\label{app:dressing}
For each $\hat\tau \in \hmf g$, the \emph{left} (resp. \emph{right})
\emph{dressing vector field} on $G$ is 
$$\XXC^l_{\hat\tau} = \pi(\theta^l_{\hat\tau}) \text{ (resp. } \XXC^r_{\hat
\tau} = -\pi(\theta^r_{\hat \tau})\text{), }$$
and $\theta^\cdot_{\hat\tau}$ is the left or right invariant $1$-form on $G$
determined by $\hat \tau$. Integrating $\XXC^\cdot_{\hat\tau}$ gives rise to a
local (global if the dressing vector fields are complete) \emph{left} (or
\emph{right}) \emph{dressing action} of \emph{the} dual group $\hat G$ on $G$,
and we say that this left (or right) dressing action consists of \emph{left} (or
\emph{right}) \emph{dressing transformations}. The Poisson Lie group $(G, \pi)$
is \emph{complete} if each left (or equivalently, right) dressing vector field
is complete. Analogously, we may define the corresponding concepts on $\hat G$.
\end{defn}
The dressing actions as defined above are the same as those in \cite{Lu, Lu1}.
Following \cite{Lu}:
\begin{defn}\label{app:poissonact}
A left action $\sigma_l: G \times P \to P$ of Poisson Lie group $(G, \pi_G)$ on
a Poisson manifold $(P, \pi_P)$ is \emph{Poisson} if $\sigma_l$ is a Poisson
map, where $G\times P$ is endowed with the product Poisson structure. Similarly
a right action $\sigma_r : P \times G \to P$ is \emph{Poisson} when $\sigma_r$
is Poisson.
\end{defn}
\begin{defn}\label{app:poissonmom} A $C^\infty$ map $\mu : P \to \hat G$ is
called a \emph{momentum mapping} for the left (resp. right) Poisson action
$\sigma_\cdot : G\times P \to P$ if for each $\tau \in \mf g = \hmf g^*$, the
infinitesimal action $X_\tau^l$ (resp. $_\tau^r$) of $\tau$ is given by
$$X_\tau^l = \pi_P(\mu^*\hat\theta^l_\tau) \text{ (resp. } X_\tau^r = -
\pi_P(\mu^*\hat\theta^r_\tau) \text{)},$$
where $\hat\theta_\tau^\cdot$ is the left (or right) invariant $1$-form on $\hat
G$ determined by $\tau$. The moment map $\mu$ of the Poisson action
$\sigma_\cdot$ is \emph{$G$-equivariant} if it's equivariant with respect to the
left (or right) dressing action of $G$ on $\hat G$.
\end{defn}
In particular, when the moment map is equivariant, we have $\mu_*(X^\cdot_\tau)
= \hat\XXC^\cdot_\tau$, where $\hat \XXC^\cdot_\tau$ is the dressing vector
field on $\hat G$ defined by $\tau \in \mf g$. Then for connected complete
Poisson Lie group $G$ (see \cite{Lu} theorem $4.8$):
\begin{theorem}\label{app:momequi}
A momentum mapping $\mu: P \to \hat G$ for a Poisson action $\sigma$ is
$G$-equivariant iff $\mu$ is a Poisson map.
\qed
\end{theorem}

\end{document}